\documentclass[reqno,11pt]{amsart}

\usepackage{amsmath}
\usepackage{amsthm}
\usepackage{mathrsfs}
\usepackage{amsfonts}
\usepackage{parskip} 
\usepackage{hyperref}
\usepackage{epsfig}

\usepackage{amsmath,amssymb,amsfonts,amsthm,amstext,verbatim}
\usepackage{enumerate,color,graphicx,
}
\usepackage{psfrag}  
\usepackage{url}
\usepackage{mdwlist}   
\usepackage{mathrsfs}
\usepackage{hyperref}
\usepackage{epsfig}
\usepackage{etoolbox}

\makeatletter
\patchcmd{\@part}{\Huge}{\normalsize}{}{} 
\makeatother

\numberwithin{equation}{section}

\allowdisplaybreaks

\newtheorem{theorem}{Theorem}[section]

\newtheorem{proposition}[theorem]{Proposition}
\newtheorem{lemma}[theorem]{Lemma}

\newtheorem{definition}[theorem]{Definition}
\newtheorem{corollary}[theorem]{Corollary}

\newcommand{\ep}{\varepsilon}

\newcommand{\R}{\mathbb{R}}
\newcommand{\N}{\mathbb{N}}
\newcommand{\Z}{\mathbb{Z}}

\newcommand{\dist}{\vert\vert}

\newcommand{\Cgjthree}{C_4}

\newcommand{\eps}{\ep}

\newcommand{\psap}{\mathsf{P}}
\newcommand{\psapleft}[1]{\mathsf{P}_{#1}^{\mathsf{left}}}

\newcommand{\psaw}{\mathsf{W}}

\newcommand{\esaw}{{\mathbb{E}_{\mathsf{W}}}}

\newcommand{\sapell}{\mathrm{SAP}^{\mathrm{left}}}
\newcommand{\sapellward}{\mathrm{SAP}^{\mathfrak{l}}}
\newcommand{\sapr}{\mathrm{SAP}^{\mathrm{right}}}
\newcommand{\saprward}{\mathrm{SAP}^{\mathfrak{r}}}

\newcommand{\fpart}{\mathrm{First}}

\newcommand{\lattv}{\Z^d}

\newcommand{\northeast}{\mathrm{NE}}
\newcommand{\eastnorth}{\mathrm{EN}}
\newcommand{\eastsouth}{\mathrm{ES}}
\newcommand{\westnorth}{\mathrm{WN}}
\newcommand{\southeast}{\mathrm{SE}}
\newcommand{\base}{\mathrm{SE}}

\newcommand{\alphamac}{\alpha'}
\newcommand{\deltamac}{\delta'}

\newcommand{\kmac}{k^*}

\newcommand\hPhi{\mathrm{High}\fpart}

\theoremstyle{remark}

\newcounter{mycount}

\def\mik{1}
\newcommand\cpsfrag[2]{\ifnum\mik=1\psfrag{#1}{#2}\fi}

\newcommand{\saw}{\mathrm{SAW}}

\newcommand{\sawfree}{\mathrm{SAW}^{*}}
\newcommand{\sap}{\mathrm{SAP}}

\newcommand{\dubbub}{\mathrm{RGJ}}

\newcommand{\maceta}{\zeta}

\newcommand{\Ch}{C_6}

\newcommand{\Cgjone}{C_2}
\newcommand{\Cgjtwo}{C_3}

\newcommand{\conindex}{C_8}

\newcommand{\epn}{\eps}
\newcommand{\epnmac}{\eps}
\newcommand{\eroom}{\epsilon_1}
\newcommand{\eclo}{\epsilon_2}
\newcommand{\esm}{\eta}

\newcommand{\thet}{\theta}

\newcommand{\hcb}{\highclose_{1/2 + \chi}}

\newcommand{\reg}{\mathsf{R}}
\newcommand{\regul}{\mathsf{K}}
\newcommand{\irr}{\mathsf{S}}
\newcommand{\regnew}[3]{\mathsf{R}_{#1,#2}^{#3}}
\newcommand{\indexc}{L}
\newcommand{\indexcj}{j}

\newcommand{\litnum}{\tfrac{1}{10}}

\newcommand{\closes}{\text{ closes}}

\newcommand\Ga{\Gamma}     
\newcommand\ga{\gamma}     

\newcommand{\hidden}[1]{}

\newcommand\strongjoin{\mathsf{GlobalMJ}}

\newcommand\globaljoin{\mathsf{GJ}}
\newcommand\ymax{y_{{\rm max}}}
\newcommand\ymin{y_{{\rm min}}}
\newcommand\xmax{x_{{\rm max}}}
\newcommand\xmin{x_{{\rm min}}}

\newcommand\righttip{\eastsouth}

\newcommand\phicl[3]{\hPhi_{#1,#2}^{#3}}

\newcommand\reverse{\overleftarrow}

\newcommand{\modify}[1]{#1_{{\rm mod}}}

\newcommand{\highclose}{\mathsf{HCP}}
\newcommand{\highpolynum}{\mathsf{HPN}}

\newcommand{\macess}{\mathsf{High}\Theta}

\newcommand{\lset}{K'}
\newcommand{\tilden}{m'}

\title[Self-avoiding polygons and walks]{On self-avoiding polygons and walks: \\ the snake method via polygon joining}
\date{}

\author[A.~Hammond]{Alan Hammond}
\address{Departments of Mathematics and Statistics, 
  U.C. Berkeley,
  Berkeley, CA, 94720-3840, U.S.A.}
\email{alanmh@stat.berkeley.edu}
\thanks{ The author is supported by NSF grant DMS-$1512908$.}

\begin{document}

\maketitle

\vspace{-0.4cm}

\begin{abstract}
For $d \geq 2$ and $n \in \N$, let $\psaw_n$ denote the uniform law on self-avoiding walks beginning at the origin in the integer lattice $\Z^d$, and write $\Gamma$ for a $\psaw_n$-distributed walk. 
We show that the {\em closing probability} $\psaw_n \big( \dist \Gamma_n \dist = 1 \big)$ that $\Gamma$'s endpoint neighbours the origin is at most $n^{-4/7 + o(1)}$ for a positive density set of odd $n$ in dimension $d = 2$.
This result is proved using the snake method, a general technique for proving closing probability upper bounds, which originated in~\cite{ontheprob} and was made explicit in~\cite{snakemethodpattern}. Our conclusion is reached by applying the snake method in unison with a polygon joining technique whose use was initiated by Madras in~\cite{Madras95}.  
\end{abstract}

\vspace{1.4cm}

\tableofcontents


\section{Introduction}

Self-avoiding walk was introduced in the 1940s by Flory and Orr \cite{Flory,Orr47} as a model of a long polymer chain in a system of such chains at very low concentration. It is well known among the basic models of discrete statistical mechanics for posing problems that are simple to state but difficult to solve. Two recent surveys are the lecture notes~\cite{BDGS11}  and \cite[Section 3]{Lawler13}. 

\subsection{The model} 

We will denote by $\N$ the set of non-negative integers.
Let $d \geq 2$. For $u \in \R^d$, let $\dist u \dist$ denote the Euclidean norm of $u$. 
A {\em walk} of length $n  \in \N$ with $n > 0$ is a map $\gamma:\{0,\cdots,n \} \to \lattv$ 
such that $\dist \gamma(i + 1) - \gamma(i) \dist = 1$ for each $i \in \{0,\cdots,n-1\}$.
An injective walk is called {\em self-avoiding}. 
A self-avoiding walk $\gamma$ of length~$n \geq 2$ is said to close (and to be closing)
if $\dist \gamma(n) \dist = 1$. 
When the {\em missing edge} connecting $\gamma(n)$ and $\gamma(0)$
is added, a polygon results.
\begin{definition}
For $n \geq 4$ an even integer, let $\ga:\{0,\ldots,n-1\} \to \Z^d$ be a closing self-avoiding walk. For  $1 \leq i \leq n-1$, let $u_i$ denote the unordered nearest neighbour edge in $\Z^d$ with endpoints $\ga(i-1)$ and $\ga(i)$. Let $u_n$ denote $\ga$'s missing edge, with endpoints $\ga(n-1)$ and~$\ga(0)$. (Note that we have excluded the case~$n=2$ so that $u_n$ is indeed not among the other $u_i$.)
We call the collection of edges 
$\big\{ u_i: 1 \leq i \leq n \big\}$ the polygon of~$\ga$. A self-avoiding polygon in $\Z^d$ is defined to be any polygon of a closing self-avoiding walk in~$\Z^d$. The polygon's length is its cardinality. 
\end{definition}
We will usually omit the adjective self-avoiding in referring to walks and polygons. Recursive and algebraic structure has been used to analyse polygons in such domains as strips, as~\cite{BMB09} describes.

Note that the polygon of a closing walk has length that exceeds the walk's by one. Polygons have even length and closing walks, odd.

Let $\saw_n$ denote the set of self-avoiding walks~$\gamma$ of length~$n$ that start at $0$, i.e., with $\gamma(0) = 0$. 
We denote by $\psaw_n$ the uniform law on $\saw_n$. 
The walk under the law $\psaw_n$ will be denoted by $\Ga$. 
The {\em closing probability} is  $\psaw_n \big( \Gamma \closes \big)$.

\subsection{A combinatorial view of the closing probability}

Let the walk number $c_n$ equal the cardinality of $\saw_n$. By equation~(1.2.10) 
of~\cite{MS93}, the limit 
$\lim_{n \in \N} c_n^{1/n}$ exists and is positive and finite; it is called the connective constant and denoted by~$\mu$, and we have~$c_n \geq \mu^n$.

Define the polygon number $p_n$ to be the number of length~$n$ polygons up to translation. By~(3.2.5) and~(3.2.9) of~\cite{MS93}, 
$\lim_{n \in 2\N} p_n^{1/n} \in (0,\infty)$ exists and coincides with~$\mu$, and $p_n \leq (d-1) \mu^n$.

The closing probability may be written in terms of the polygon and walk numbers.
 There are $2n$ closing walks whose polygon is a given polygon of length~$n$, since there are $n$ choices of missing edge and two of orientation. Thus,
\begin{equation}\label{e.closepc}
\psaw_n \big( \Gamma \closes \big) = \frac{2(n+1) p_{n+1}}{c_n} \, ,
\end{equation}
for any $n \in \N$ (but non-trivially only for odd values of $n$).
 
\subsection{Corrections to exponential growth} 
 
We
define the real-valued {\em polygon number deficit} and {\em walk number excess} exponents $\thet_n$ and~$\xi_n$ according to the formulas
\begin{equation}\label{e.thetn}
 p_n = n^{-\thet_n} \cdot \mu^n \, \, \, \textrm{for $n \in 2\N$} 
\end{equation}
and
\begin{equation}\label{e.xin}
 c_n = n^{\xi_n} \cdot \mu^n \, \, \, \textrm{for $n \in \N$} \, .
\end{equation}
Since $c_n \geq \mu^n$ and $p_n \leq (d-1) \mu^n$,  $\liminf_{n \in 2\N} \thet_n$ and $\inf_{n \in \N} \xi_n$ are non-negative in any dimension~$d \geq 2$.

\subsection{Closing probability upper bounds: a review}

An interesting first task regarding closing probability upper bounds is to show that $\psaw_n \big( \Gamma_n \closes \big) \to 0$
as $n \to \infty$ for any $d \geq 2$; this question was posed by Itai Benjamini to the authors of~\cite{ontheprob} who showed an upper bound of $n^{-1/4 - o(1)}$ whenever $d \geq 2$. One way of approaching the question -- not the method of~\cite{ontheprob} -- is to improve the above stated $\liminf_{n \in 2\N} \big( \xi_{n-1} + \theta_n \big) \geq 0$ so that the right-hand side takes the form $1 + \chi$, for some $\chi > 0$, in which case, a closing probability upper bound of the form $n^{-\chi + o(1)}$ will result from~(\ref{e.closepc}).

 Hara and Slade~\cite[Theorem 1.1]{HS92} used the lace expansion to prove that $c_n \sim C \mu^n$ when $d \geq 5$, but the problem of showing that $c_n \mu^{-n}$ tends to infinity as $n \to \infty$ appears to be open when $d$ equals two or three (in which cases, it is anticipated to be true). As such, this combinatorial approach to closing probability upper bounds has depended on advances on the polygon rather than the walk side. Madras and Slade \cite[Theorem 6.1.3]{MS93} have proved that $\thet_n \geq d/2 + 1$ when $d \geq 5$ for {\em spread-out} models, in which the vertices of $\Z^d$ are connected by edges below some bounded distance. Thus, the conclusion that the closing probability decays as fast as $n^{-d/2}$ has been reached for such models when $d \geq 5$. (In passing: this conclusion may be expected to be sharp, but the opposing lower bound is not known to the best of the author's knowledge.) When the model is nearest neighbour, the upper bound is in essence known for $d \geq 5$ in an averaged sense:~\cite[Theorem~$1.3$]{HS92} states that $\sum_n n^a \psaw_n (\Ga \closes) < \infty$ for any $a < d/2 - 1$; in fact, the counterpart to this closing probability bound  for a walk ending at any given displacement from the origin is also proved, uniformly in this displacement.

It is interesting to note that, while the problem of showing that $\limsup_{n \in 2\N} \theta_n$ is finite is open, it is known from classical results that, for any dimension $d \geq 2$, 
$\liminf_{n \in 2\N} \theta_n < \infty$. Indeed, a technique~\cite[Theorem~$3.2.4$]{MS93} that rearranges two self-avoiding bridges into a closing walk (and, by the addition of the missing edge, into a polygon) combines with the subsequential availability of bridges of length~$n$ 
(meaning $b_n \mu^{-n} \geq n^{-1-o(1)}$ for infinitely many $n$), an availability which is   due to the divergence of the bridge generating function at its critical point~\cite[Corollary~$3.1.8$]{MS93}, to yield the inference that $\liminf_{n \in 2\N} \theta_n \leq d+5$.  

In dimension $d=2$, Madras~\cite{Madras95} showed that $\liminf \theta_n \geq 1/2$
in 1995. To show this, he introduced a polygon joining technique 
of which we will also make use in the present article (and which is reviewed in Section~\ref{s.madrasjoin}). His technique has recently been exploited to prove the next result,~\cite[Theorem~$1.3$]{ThetaBound}.
\begin{definition}
The limit supremum density of 
a set $A$ of even, or odd, integers~is
$$
\limsup_n \frac{\big\vert A \cap [0,n] \big\vert}{\vert 2\N \cap [0,n] \vert} = \limsup_n \, n^{-1} \big\vert A \cap [0,2n] \big\vert \, .
$$
When the corresponding limit infimum density equals the limit supremum density, we naturally call it the limiting density. 
\end{definition}
\begin{theorem}\label{t.polydev}
 Let $d=2$. For any $\delta > 0$,
the  limiting density of the set of $n \in 2\N$ for which $\thet_n \geq 3/2 - \delta$ is equal to one. 
\end{theorem}

As such, the combinatorial perspective offered by~(\ref{e.closepc}) alongside polygon joining via Theorem~\ref{t.polydev}
has yielded the inference that $\psaw_n \big( \Gamma \closes \big) \leq n^{-1/2 + o(1)}$ on a full density set of odd~$n$ when $d=2$. 

This bound complements the $n^{-1/4 + o(1)}$ general dimensional upper bound from~\cite{ontheprob}. In fact, the method of~\cite{ontheprob} was reworked in~\cite{snakemethodpattern}
to achieve the following  $n^{-1/2 + o(1)}$ bound.   
\begin{theorem}\label{t.closingprob}
Let $d \geq 2$. For any $\eps > 0$ and $n \in 2\N + 1$ sufficiently high,
$$
\psaw_n \big( \Gamma \closes \big) \leq n^{-1/2 + \eps} \, .
$$
\end{theorem}
The method of proof of this result was not a combinatorial analysis via~(\ref{e.closepc}) but rather a probabilistic approach called the snake method (about which more shortly).

We may summarise this review of existing closing probability upper bounds as follows. An upper bound of the form $n^{-1/2 + o(1)}$ in dimension~$d=2$ has been obtained in two very different ways. 
Via combinatorics and polygon joining, Theorem~\ref{t.polydev}
asserts this bound for typical odd~$n$. Via the probabilistic snake method, Theorem~\ref{t.closingprob} asserts the bound for each  $d \geq 2$ and for all high $n$.

\subsection{The main result}

This article's principal conclusion strengthens this upper bound when $d=2$ beyond the threshold of one-half in the two existing methods.

\begin{theorem}\label{t.thetexist}
Let $d = 2$. 
\begin{enumerate} 
\item
For any $\eps > 0$, the bound 
$$
\psaw_n \big( \Gamma \closes \big) \leq n^{-4/7 + \eps} 
$$
holds on a set of $n \in 2\N + 1$ of limit supremum density at least $1/{1250}$. 
\item
Suppose that the limits $\thet : = \lim_{n \in 2\N} \thet_n$
and $\xi : = \lim_{n \in \N} \xi_n$
exist in $[0,\infty]$. 
Then $\thet + \xi \geq 5/3$. Since 
\begin{equation}\label{e.closepc.formula}
\psaw_n \big( \Gamma \closes \big) = n^{-\theta - \xi + 1 + o(1)} 
\end{equation}
as $n \to \infty$ through odd values of $n$
by (\ref{e.closepc}), the closing probability is seen to be bounded above by $n^{-2/3 + o(1)}$. 
\end{enumerate} 
\end{theorem} 
(When $\thet + \xi = \infty$, (\ref{e.closepc.formula}) should be interpreted as asserting a superpolynomial decay in $n$ for the left-hand side.)

However far from rigorous the hypotheses of Theorem~\ref{t.thetexist}(2)  may be, their validity is uncontroversial. For example, the limiting value $\thet$ is predicted to exist and to satisfy a relation with the Flory exponent~$\nu$ for mean-squared radius of gyration. The latter exponent is specified by the putative formula $\esaw_n \, \dist \Gamma(n) \dist^2   = n^{2\nu + o(1)}$,
where $\esaw_n$ denotes the expectation associated with~$\psaw_n$ (and where note that $\Gamma(n)$ is the non-origin endpoint of $\Gamma$); in essence, $\dist \Ga(n) \dist$  is supposed to be typically of order $n^{\nu}$. 
The hyperscaling relation that is expected to hold between $\thet$ and $\nu$ is $\thet   =   1 +  d \nu$
where the dimension $d \geq 2$ is arbitrary. In $d = 2$, $\nu = 3/4$ and thus $\theta = 5/2$ is expected. That $\nu = 3/4$ was predicted by the Coulomb gas formalism \cite{Nie82,Nie84} and then by conformal field theory \cite{Dup89,Dup90}. 
We mention also that $\xi = 11/32$ is expected when $d=2$; in light of the $\thet = 5/2$ prediction and~(\ref{e.closepc}),
 $\psaw_n \big( \Ga \closes\big) = n^{-\psi + o(1)}$ with $\psi = 59/32$ is expected. The $11/32$ value was predicted by Nienhuis in \cite{Nie82} and can also be deduced from calculations concerning SLE$_{8/3}$: see \cite[Prediction~5]{LSW}. We will refer to $\psi$ as the {\em closing} exponent.

\subsection{The snake method and  polygon joining in unison}

We have mentioned that the snake method is a probabilistic tool for proving closing probability upper bounds.  It is introduced in~\cite{snakemethodpattern}; (it originated, though was not made explicit, in~\cite{ontheprob}). The method states a condition that must be verified in order to achieve some such upper bound. In~\cite{snakemethodpattern}, this condition was verified with the use of a technique we may call Gaussian pattern fluctuation, in order to obtain Theorem~\ref{t.closingprob}.

This article is devoted to proving Theorem~\ref{t.thetexist}.
We will obtain the result by  a second use of the snake method. In order to verify the necessary condition, we will replace the use of pattern fluctuation with the technique of polygon joining seen in \cite{Madras95,ThetaBound}.

That is, Theorem~\ref{t.thetexist} 
pushes above the threshold of one-half for the exponent in the closing probability upper bound by combining the combinatorial and probabilistic perspectives that had been used in the two proofs that led to the lower border of this threshold. The first part of the theorem is probably the more significant conclusion.  Theorem~\ref{t.thetexist}(2) is only a conditional result, but it serves a valuable expository purpose: its proof is that of the theorem's first part with certain technicalities absent.

\medskip

\noindent{\bf Structure of the paper.}
It will by now be apparent that  some groundwork is needed before the proof of Theorem~\ref{t.thetexist} may be given. The two elements that must be reviewed (and analysed) are the polygon joining technique and the snake method. After some general notation in Section~\ref{s.notation},  the respective reviews are contained in 
Sections~\ref{s.polyjoin},~\ref{s.madrasjoin} and~\ref{s.pjprep} and in Section~\ref{s.five}.  The proof of the main result then appears in
the remaining eight sections, of which the first, Section~\ref{s.expo}, offers an exposition. 
See the beginning of Section~\ref{s.polyjoin}
for a summary of the content of the joining material in the ensuing three sections and the beginning of  Section~\ref{s.theoremsecond} for a summary of the structure of the sections giving the proof of Theorem~\ref{t.thetexist}.

\medskip

\noindent{\bf A few words on possible reading of related works.}
This article is intended to be read on its own.
It is nonetheless the case that, since this article makes use of  a fairly subtle combination of existing techniques, the reader may benefit from reading some of the related material. For example, the application of the snake method made here is rather more technical than the one in~\cite{snakemethodpattern}, and it is imaginable that reading parts of that 
article might bring the reader a certain familiarity with the method that would be a useful (though not necessary) preparation for reading the present one. Our polygon joining method refines that in~\cite{ThetaBound} and again some review of that work is not necessary but possibly helpful. We also mention that the online article~\cite{CJC} gathers together the content of the three articles \cite{ThetaBound},~\cite{snakemethodpattern} and~the present work, while also providing some further exposition and highlighting some other connections between the concepts in the articles; the reader who wishes to study all the various theorems collectively is encouraged to consult this work. 

Alternatively, optional further reading may be selected by noting   the outside inputs used by this paper. Theorem~\ref{t.polydev}
and Corollary~\ref{c.globaljoin} are proved in~\cite{ThetaBound}: they are 
Theorem~$1.3$ and Corollary~$4.6$ of that article.
Lemma~\ref{l.closecard} and Theorem~\ref{t.snakesecondstep} appear in~\cite{snakemethodpattern}: see Lemma~$2.2$ and Theorem~$3.2$ in that work. The latter two results concern the application and setup of the snake method. It is worth pointing out that their proofs are respectively almost trivial and fairly short (three to four pages in the latter case).
It is thus suggested to the reader who wishes to make a 
moderate attempt at understanding the paper's inputs to read the proofs of these results from Sections~$2.2$ and~$3.2$ of~\cite{snakemethodpattern} when their statements are reached. 

\medskip

\noindent{\bf Acknowledgments.} 
I am very grateful to a referee for a thorough discussion of an earlier version of~\cite{CJC}. Indeed, the present form of  Theorem~\ref{t.thetexist}(1) is possible on the basis of a suggestion made by the referee that led to a strengthening of the original version of~\cite[Theorem 1.3]{ThetaBound}. I thank a second referee for a thorough reading and valuable comments.  
I thank Hugo Duminil-Copin and Ioan Manolescu for many stimulating and valuable conversations about the central ideas in the paper. I thank Wenpin Tang for useful comments on a draft version of~\cite{CJC}. I would also like to thank Itai Benjamini for suggesting the problem of studying upper bounds on the closing probability.

\section{Some general notation and tools}\label{s.notation}

\subsubsection{Cardinality of a finite set $A$} This is denoted by either $\# A$ or $\vert A \vert$.

\subsubsection{Multi-valued maps}\label{s.mvp}
For a finite set $B$, let $\mathcal{P}(B)$ denote its power set.
Let $A$ be another finite set. A {\em multi-valued map} from $A$ to $B$ is a function $\Psi: A \to \mathcal{P}(B)$. An arrow is a pair $(a,b) \in A \times B$ for which $b \in \Psi(a)$; such an arrow is said to be outgoing from $a$ and incoming to $b$. We consider multi-valued maps in order to find lower bounds on $\vert B \vert$, and for this, we need upper (and lower) bounds on the number of incoming (and outgoing) arrows. This next lemma is an example of such a lower bound.
\begin{lemma}\label{l.mvm}
Let $\Psi: A \to \mathcal{P}(B)$.
Set $m$ to be the minimum over $a \in A$ of the number of arrows outgoing from $a$, and $M$ to be the maximum over $b \in B$ of the number of arrows incoming to~$b$. Then
  $\vert B \vert \geq m M^{-1} \vert 
A \vert$.
\end{lemma}
\noindent{\bf Proof.} The quantities $M \vert B \vert$ and $m \vert A \vert$
are upper and lower bounds on the total number of arrows. \qed

\subsubsection{Denoting walk vertices and subpaths}
For $i,j \in \N$ with $i \leq j$, we write $[i,j]$ for  $\big\{ k \in \N:  i \leq k \leq j \big\}$. For a walk $\ga:[0,n] \to \Z^d$ and $j \in [0,n]$, we write $\ga_j$ in place of $\ga(j)$. For $0 \leq i \leq j \leq n$, $\ga_{[i,j]}$ denotes the subpath $\ga_{[i,j]}: [i,j] \to \Z^d$ given by restricting $\ga$. 

\subsubsection{Notation for certain corners of polygons}

\begin{definition}\label{d.corners}
The Cartesian unit vectors are denoted by $e_1$ and $e_2$ and the coordinates of $u \in \Z^2$ by $x(u)$ and $y(u)$.
For a finite set of vertices $V \subseteq \Z^2$, we define the northeast vertex $\northeast(V)$ in $V$ to be that element of $V$ of maximal $e_2$-coordinate; should there be several such elements, we take $\northeast(V)$ to be the one of maximal $e_1$-coordinate. That is, $\northeast(V)$ is the uppermost element of $V$, and the rightmost among such uppermost elements if there are more than one. (Rightmost means of maximal $e_1$-coordinate.) Using the four compass directions, we may similarly define eight elements of $V$, 
including the lexicographically minimal and maximal elements of $V$, $\mathsf{WS}(V)$ and $\mathsf{EN}(V)$. We extend the notation to any self-avoiding walk or polygon $\gamma$, writing for example $\northeast(\gamma)$ for $\northeast(V)$, where $V$ is the vertex set of $\gamma$. 
For a polygon or walk $\ga$, set $\ymax(\ga) = y\big(\northeast(\gamma)\big)$, $\ymin(\ga) = y\big(\southeast(\gamma)\big)$, $\xmax(\ga) = x\big(\eastnorth(\gamma)\big)$ and $\xmin(\ga) = x\big(\westnorth(\gamma)\big)$.
The height~$h(\ga)$ of $\ga$ 
 is $\ymax(\gamma) - \ymin(\gamma)$  and its width~$w(\gamma)$ is  $\xmax(\gamma) - \xmin(\gamma)$.
\end{definition}

\subsubsection{Polygons with northeast vertex at the origin}

For $n \in 2\N$, let $\sap_n$ denote the set of length $n$ polygons $\phi$ such that $\northeast(\phi) = 0$. The set $\sap_n$ is in bijection with equivalence classes of length~$n$ polygons where polygons are identified if one is a translate of the other. Thus, $p_n =   \vert \sap_n \vert$. 

We write $\psap_n$ for the uniform law on $\sap_n$. 
A polygon  sampled with law $\psap_n$ will be denoted by $\Ga$, as a walk with law $\psaw_n$ is.

There are $2n$ ways of tracing the vertex set of a polygon $\phi$ of length~$n$: $n$ choices of starting point and two of orientation. We now select one of these ways. Abusing notation, we may write $\phi$ as a map from $[0,n]$ to $\Z^2$, setting $\phi_0 = \northeast(\phi)$, $\phi_1 = \northeast(\phi) - e_1$, and successively defining $\phi_j$ to be the previously unselected vertex for which $\phi_{j-1}$ and $\phi_j$ form the vertices incident to an edge in $\phi$, with the final choice $\phi_n = \northeast(\phi)$ being made. Note that $\phi_{n-1} = \northeast(\phi) - e_2$.


\subsubsection{Plaquettes}

The shortest non-empty polygons contain four edges. Certain such polygons play an important role in several arguments and we introduce notation for them now.

\begin{definition}\label{d.plaquette}
A plaquette is a polygon with four edges. Let $\phi$ be a polygon. A plaquette $P$ is called a join plaquette of~$\phi$ if $\phi$ and $P$ intersect at precisely the two horizontal edges of $P$. 
Note that when $P$ is a join plaquette of $\phi$, the operation of removing the two horizontal edges in $P$ from $\phi$ and then adding in the two vertical edges in $P$ to~$\phi$ results in two disjoint polygons whose lengths sum to the length of $\phi$. We use symmetric difference notation and denote the output of this operation by $\phi \, \Delta \, P$.

The operation may also be applied in reverse: for two disjoint polygons $\phi^1$ and $\phi^2$, each of which contains one vertical edge of a plaquette~$P$, the outcome $\big( \phi^1 \cup \phi^2 \big) \, \Delta \, P$ of removing $P$'s vertical edges and adding in 
its horizontal ones is a polygon whose length is the sum of the lengths of $\phi^1$ and $\phi^2$.
\end{definition}

\section{A very brief overview of polygon joining}\label{s.polyjoin}

The next three sections define and explain our method of polygon joining. In the present section, we set the scene with a brief review of pertinent concepts. Section~\ref{s.madrasjoin} reviews the local details of Madras' joining procedure. Section~\ref{s.pjprep}
introduces and analyses the detailed mechanism for joining that we will use to prove Theorem~\ref{t.thetexist}.

Recall that the two-dimensional Theorem~\ref{t.polydev} is proved in~\cite{ThetaBound} via polygon joining. 
Since it may be conceptually useful, we now recall some elements of~\cite{ThetaBound} and give a very brief overview of the method of proof of its conclusion that $\theta_n \geq 3/2 - o(1)$ for typical~$n$. Section 3 of~\cite{ThetaBound} is a heuristical discussion of polygon joining which expands on the next few paragraphs. We review three joining arguments, the third one leading to the above conclusion. 

Consider first a pair of polygons in $d=2$ of say equal length $n$. The second may be pushed to a location where it is disjoint from, and to the right of, the first, so that there exists a plaquette whose left vertical edge lies in the first polygon and whose right vertical edge lies in the second. The reverse operation 
in the latter paragraph of Definition~\ref{d.plaquette} may be applied to this plaquette and a single polygon of length $2n$ results. A moment's thought shows that, since the original pair of polygons had equal length, the location of the concerned plaquette is identifiable from the image polygon. Thus, polygon joining in its simplest guise yields the polygon superadditivity bound $p_{2n} \geq p_n^2$. 

The next polygon joining argument is due to Madras~\cite{Madras95}. 
The first polygon may be oriented to have height at least $n^{1/2}$. Vertically displacing the second before attempting joining should yield order $n^{1/2}$ locations for joining and thus $p_{2n} \geq n^{1/2} p_n^2$ and in essence $\theta_n \geq 1/2$.

The third joining argument is performed in \cite{ThetaBound}: working with $n/2$ in place of $n$, we vary the length pair for the concerned polygons from $(n/2,n/2)$ to $(n/2-j,n/2+j)$ where $j$ varies over choices such that $\vert j \vert \leq n/4$, and we seem to obtain 
\begin{equation}\label{e.polyjoin}
p_n \geq n^{1/2 - o(1)} \sum_{j = -n/4}^{n/4} p_{n/2-j} p_{n/2+j}
\end{equation}
 whence $\liminf_{m \in 2\N} \theta_m \geq 3/2 - o(1)$. In fact, some care is needed to derive this bound. To obtain it, it must be the case that a typical length-$n$ polygon formed by such joinings is such that the plaquette used to join is identifiable among at most~$n^{o(1)}$ possible locations (we can call this circumstance {\em efficient plaquette identification}). The join plaquette in question is  `global' in the rough sense that it borders two macroscopic chambers.
Here is the precise definition from~\cite{ThetaBound}, of which we will again make use. 
\begin{definition}\label{d.globaljoin}
For $n \in 2\N$, let $\phi \in \sap_n$. 
A join plaquette $P$ of $\phi$ is called {\em global} if the two polygons 
comprising $\phi \, \Delta \, P$ may be labelled $\phi^\ell$ and $\phi^r$ in such a way that 
\begin{itemize}
\item every rightmost vertex in $\phi$ is a vertex of $\phi^r$;
\item and  $\northeast(\phi)$ is a vertex of $\phi^\ell$.
\end{itemize}
Write $\globaljoin_\phi$ for the set of global join plaquettes of the 
polygon~$\phi$.
\end{definition}

\begin{figure}
    \begin{center}
      \includegraphics[width=0.75\textwidth]{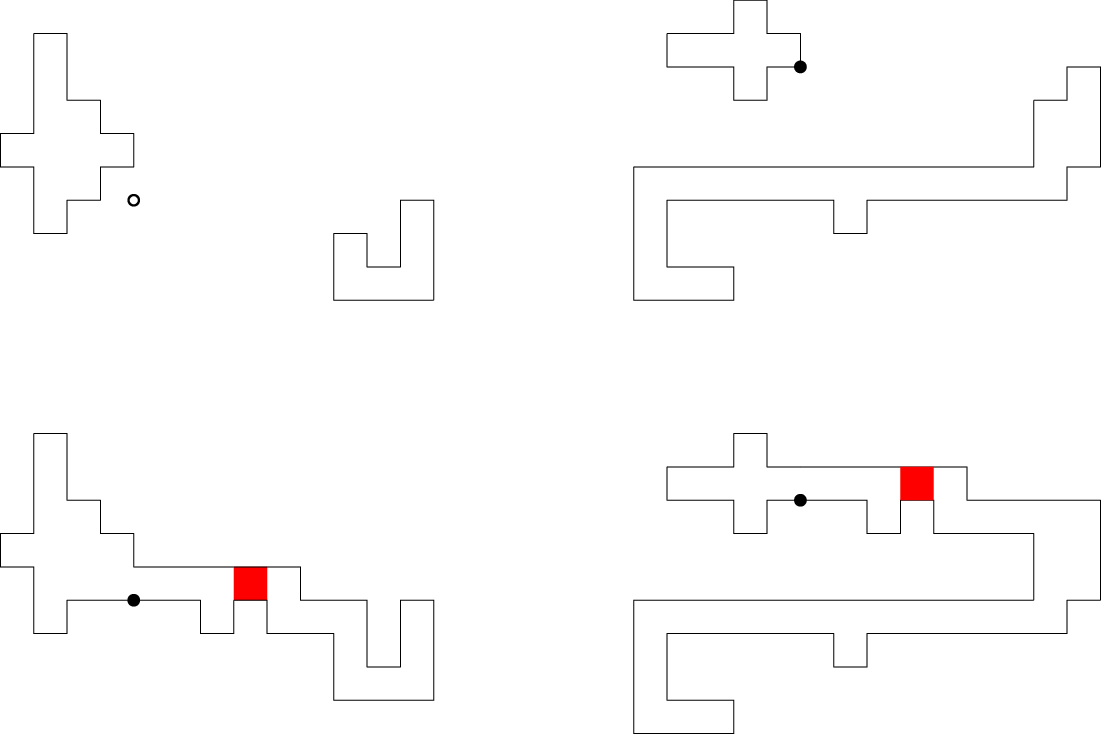}
    \end{center}
    \caption{Two pairs of Madras joinable polygons and the Madras join polygon of each pair. The junction plaquette is shaded red (or, a touch pedantically, the interior of the union of its constituent edges is), while the circle or disk indicates an aspect of the detailed mechanism of Madras joining omitted in our overview.}\label{f.madrasjoinable}
  \end{figure}

For example, the four edges that form the boundary of the unit square shaded red in each of the polygons~$\phi$ in the lower part of Figure~\ref{f.madrasjoinable} are a plaquette belonging to $\globaljoin_\phi$. 

In~\cite{ThetaBound}, efficient plaquette identification is proved to work only in the case that the image polygon length lies in a set $\highpolynum$ of indices of {\em high polygon number}. 
Next we state Definition~4.1 and Corollary~4.6 of~\cite{ThetaBound}.

\begin{definition}\label{d.highpolynum}
For $\maceta > 0$, the set  $\highpolynum_\maceta \subseteq 2\N$ of $\maceta$-{\em high polygon number} indices is given by
$$
 \highpolynum_\maceta = \Big\{ n \in 2\N : p_n \geq n^{- \maceta} \mu^n \Big\} \, .
$$
\end{definition}

\begin{corollary}\label{c.globaljoin}
For all $\maceta > 0$ and $\Cgjone > 0$, there exist $\Cgjtwo, \Cgjthree > 0$ such that, for $n \in \highpolynum_\maceta$,
$$
\psap_{n} \Big( \, \big\vert \globaljoin_\Gamma \big\vert \geq \Cgjtwo \log n \Big)
 \leq \Cgjthree n^{- \Cgjone} \, .
$$
\end{corollary}
Thus, when $n \in \highpolynum_\maceta$, a form of~(\ref{e.polyjoin}) is obtained, in \cite[Proposition 4.2]{ThetaBound}.

Many of these ideas will also be needed in the present article, which is why we have just informally recounted them.

\section{Madras' polygon joining procedure}\label{s.madrasjoin}

When a pair of polygons is close, there may not be a plaquette whose vertical edges are divided between the two elements of the pair in the manner discussed above. Local details of the structure of the two polygons in the locale of near contact must be broken into cases and a suitable local modification made in each one, in order to join the pair of polygons.  Madras~\cite{Madras95} defined such a joining technique which works in a general way. 
For the most part, the precise details of this mechanism need not be the focus of the reader's attention; it is only in the proof of one technical result, Lemma~\ref{l.leftlongbasic}, that they must be analysed, (and even then the result is plausible without this analysis). In this section, we explain in rough terms the principal relevant features. The reader who would like to see a detailed overview is invited to read \cite[Section 4.1]{ThetaBound} where Madras' polygon joining technique is reviewed in detail and the specifics of the notation that we use in relation to it are given; a figure depicting the various cases also appears there.

In outline, then, here is Madras' technique:
consider two polygons $\tau$ and $\sigma$ of lengths $n$ and $m$ for which the intervals 
$$
\textrm{$\big[ \ymin(\tau) - 1 , \ymax(\tau) + 1 \big]$ and  $\big[ \ymin(\sigma) - 1 ,\ymax(\sigma) + 1 \big]$ intersect} \, .
$$
 Madras' procedure joins  $\tau$ and $\sigma$ to form a new polygon of length $n + m + 16$ in the following manner. The polygon $\sigma$ is first pushed far to the right and then back to the left until it is a few units from touching $\tau$. A local modification is made to both polygons around the contact point, with a net gain in length to both polygons of eight. The modification causes the polygons to reach closer to one another horizontally by four, five or six units and a {\em horizontal translational adjustment} of the right polygon in one or other direction may be needed in order to obtain the desired outcome, which is a polygon pair for which there exists a plaquette in the locale of modification with left vertical edge in the left polygon and right vertical edge in the right polygon. The reverse operation in Definition~\ref{d.plaquette} applied to this plaquette then yields a joined polygon of length $n + 16$. This polygon is defined to be the {\em Madras join polygon} $J(\tau,\sigma)$. The plaquette used at the final step of the construction of this polygon is called the {\em junction} plaquette.
 
Sometimes a pair of polygons are located in such a relative position (with the right polygon a few units away from the left near the point of contact) that in the above procedure, no horizontal adjustment is needed. In this case, the pair is called {\em Madras joinable}.

\section{Polygon joining in preparation for the proof of Theorem~\ref{t.thetexist}}\label{s.pjprep}

 Our purpose in this section is to specify the polygon pairs that we will join to prove Theorem~\ref{s.pjprep} and discuss some properties of the joining.

It may be apparent from the review in Section~\ref{s.polyjoin} that, in the proof of Theorem~\ref{t.polydev} that appears in \cite{ThetaBound}, certain {\em left} polygons are Madras joined to other {\em right} polygons. It serves our purpose of specifying the polygon joining details needed for Theorem~\ref{t.thetexist}'s proof to recall the definitions and some properties of these polygons, and we do so in the first of this section's six subsections.
In the second subsection, we refine the definitions of these polygon types, and also of the joining mechanism for them. The polygons that are formed under the new, stricter, mechanism will be called 
{\em regulation global join} polygons: in the third subsection, these polygons are defined and their three key properties stated.
The remaining three subsections prove these properties in turn.

\subsection{Left and right polygons}

Let $\phi$ be a polygon. Recall from Definition~\ref{d.corners} the notation $\ymax(\phi)$ and $\ymin(\phi)$, as well as the height~$h(\phi)$ and width~$w(\phi)$.

\begin{definition}\label{d.leftright}
For $n \in 2\N$, let $\sapellward_n$ denote the set of {\em left} polygons $\phi \in \sap_n$
such that
\begin{itemize}
\item $h(\phi) \geq w(\phi)$ (and thus, by a trivial argument, $h(\phi) \geq n^{1/2}$),
\item and $y\big(\righttip(\phi) \big) \leq \tfrac{1}{2} \big( \ymin(\phi) + \ymax(\phi) \big)$. 
\end{itemize}
Let $\saprward_n$ denote the set of {\em right}  polygons $\phi \in \sap_n$
such that
\begin{itemize}
\item $h(\phi) \geq w(\phi)$.
\end{itemize}
\end{definition}
\begin{lemma}\label{l.polysetbound.newwork}
For $n \in 2\N$,
$$
\big\vert \sapellward_n \big\vert \geq \tfrac{1}{4} \cdot \big\vert \sap_n \big\vert \, \, \, 
\textrm{and}  \, \, \,
\big\vert \saprward_n \big\vert \geq  \tfrac{1}{2} \cdot \big\vert \sap_n \big\vert \, .
$$
\end{lemma}
\noindent{\bf Proof.} An  element $\phi \in \sap_n$
not in $\saprward_n$ is brought into this set
by right-angled rotation. 
If, after the possible rotation, it is not in $\sapellward_n$, it may 
brought there by reflection in the $x$-axis. \qed

\begin{definition}\label{d.globallymj}
A Madras joinable polygon pair $(\phi^1,\phi^2)$
is called globally Madras joinable if the junction plaquette of the join polygon $J(\phi^1,\phi^2)$ is a global join plaquette of  $J(\phi^1,\phi^2)$.
\end{definition}
Both polygon pairs in the upper part of Figure~\ref{f.madrasjoinable} are globally Madras joinable.

The next result is \cite[Lemma 4.11]{ThetaBound}.
\begin{lemma}\label{l.strongjoin}
Let $n,m \in 2\N$ and let $\phi^1 \in \sapellward_n$ and $\phi^2 \in \saprward_m$.

Every value 
\begin{equation}\label{e.kvalues}
k \in \Big[ y \big(\righttip(\phi^1)\big)  \, , \,  y\big(\righttip(\phi^1)\big)  +  \min \big\{ n^{1/2}/2,m^{1/2} \big\}  - 1 \Big]
\end{equation}
is such that $\phi^1$ and some horizontal shift of $\phi^2 + ke_2$ is globally Madras joinable. 

Write $\strongjoin_{(\phi^1,\phi^2)}$ for the set of $\vec{u} \in \Z^2$ such that 
the pair $\phi^1$ and $\phi^2 + \vec{u}$ is globally Madras joinable.
Then
$$
  \big\vert \strongjoin_{(\phi^1,\phi^2)} \big\vert \geq   \min \big\{n^{1/2}/2,m^{1/2} \big\}  \, .
$$
\end{lemma}

  \begin{figure}
    \begin{center}
      \includegraphics[width=0.75\textwidth]{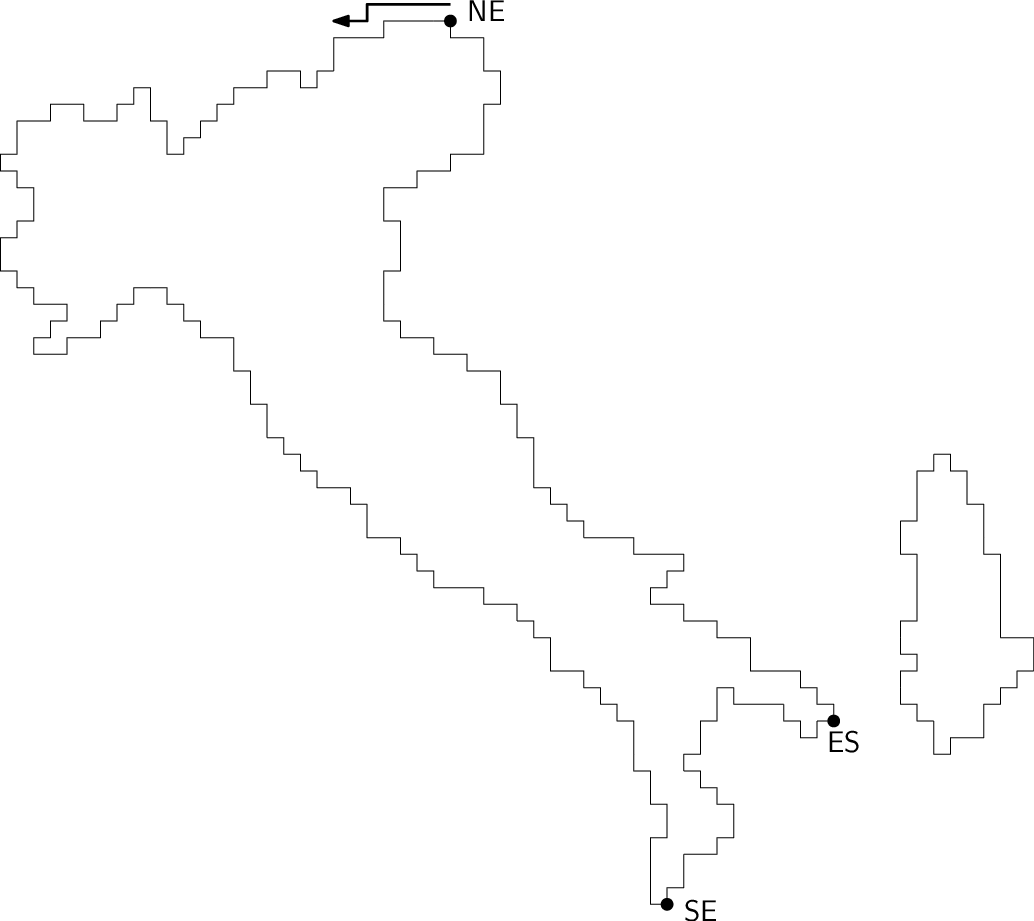}
    \end{center}
    \caption{Illustrating Lemma~\ref{l.strongjoin} with two polygons, $\mathsf{ITA}$ and $\mathsf{ALB}$. Denoting the polygons' lengths by $n$ and $m$, note that the translates of  $\mathsf{ITA}$ and $\mathsf{ALB}$ with $\northeast = 0$ belong to $\sapellward_n$ and $\saprward_m$. 
The start of the counterclockwise tour of $\mathsf{ITA}$ from $\northeast(\mathsf{ITA})$ 
dictated by convention is marked by an arrow.    
    The depicted polygons are globally Madras joinable after a horizontal shift of $\mathsf{ALB}$, with the necessary surgery to $\mathsf{ITA}$ occurring in a vicinity of $\eastsouth(\mathsf{ITA})$.
Vertical shifts of $\mathsf{ALB}$
that leave this polygon in an easterly line of sight from  $\eastsouth(\mathsf{ITA})$ will maintain this state of affairs.}
    \label{fig:albita}
  \end{figure}

\subsection{Left and right polygons revisited}\label{s.refpert}

We will refine these notions now, specifying new left and right polygon sets $\sapell_n$
and $\sapr_n$. Each is a subset of its earlier counterpart. 

In order to specify $\sapell_n$, we now define a {\em left-long} polygon.
Let $\phi \in \sap_n$. 
We employ the notationally abusive parametrization $\phi:[0,n] \to \Z^2$ such that $\phi_0 = \phi_n = \northeast(\phi)$ and $\phi_1 = \northeast(\phi) - e_1$. If $j \in [0,n]$ is such that $\phi_j = \base(\phi)$, note that $\phi$ may be partitioned into two paths $\phi_{[0,j]}$ and $\phi_{[j,n]}$. 
The two paths are edge-disjoint, and the first lies to the left of the second: indeed, writing~$H$ for the horizontal strip $\big\{ (x,y) \in \R^2: y\big(\base(\phi)\big) \leq y \leq y\big(\northeast(\phi)\big) \big\}$, the set $H \setminus \cup_{i=0}^{j-1} [\phi_i,\phi_{i+1}]$ has two  unbounded connected components; one of these, the right component, contains $H \cap \big\{ (x,y) \in \R^2: x \geq C \big\}$ for large enough $C$; and the union of the edges $[\phi_i,\phi_{i+1}]$, $j \leq i \leq n-1$, excepting the points $\phi_j = \base(\phi)$ and $\phi_n = \northeast(\phi)$, lies in the right component. It is thus natural to call $\phi_{[0,j]}$ the left path, and $\phi_{[j,n]}$, the right path. We call $\phi$ {\em left-long}  if $j \geq n/2$ and {\em right-long} if $j \leq n/2$.  

For $n \in 2\N$, let $\sapell_n$ denote the set of left-long elements of
$\sapellward_n$.
We simply set $\sapr_n$ 
equal to $\saprward_n$.
Note that $\mathsf{ITA}$ in Figure~\ref{fig:albita} is an element of $\sapell_n$.

\begin{lemma}\label{l.polysetbound}
For $n \in 2\N$,
$$
\big\vert \sapell_n \big\vert \geq \tfrac{1}{8} \cdot \big\vert \sap_n \big\vert \, \, \, 
\textrm{and}  \, \, \, 
\big\vert \sapr_n \big\vert \geq  \tfrac{1}{2} \cdot \big\vert \sap_n \big\vert \, .
$$
\end{lemma}
\noindent{\bf Proof.} The second assertion has been proved in Lemma~\ref{l.polysetbound.newwork}. In regard to the first, write the three requirements for an element $\phi \in \sap_n$ to satisfy $\phi \in \sapell_n$
in the order: $h(\phi) \geq w(\phi)$; $\phi$ is left-long; and  $y\big(\righttip(\phi) \big) \leq \tfrac{1}{2} \big( \ymin(\phi) + \ymax(\phi) \big)$.
 These conditions may be satisfied as follows.
\begin{itemize}
\item The first property may be ensured by a right-angled counterclockwise rotation if it does not already hold. 
\item It is easy to verify that, when a polygon is reflected in a vertical line, the right path is mapped to become a sub-path of the left path of the image polygon. Thus, a right-long polygon maps to a left-long polygon under such a reflection.
(We note incidentally the asymmetry in the definition of left- and right-long: this last statement is not always true {\em vice versa}.) A polygon's height and width are unchanged by either horizontal or vertical reflection, so the first property is maintained if a reflection is undertaken at this second step.
\item The third property may be ensured if necessary by reflection in a horizontal line.  The first property's occurrence is not disrupted for the reason just noted. Could the reflection disrupt the second property? 
When a polygon is reflected in a horizontal line, $\northeast$ and $\base$ in the domain  map to $\base$ and $\northeast$ in the image. The left path maps to the reversal of the left path, (and similarly for the right path). Thus, any left-long polygon remains left-long when it is reflected in such a way. We see that the second property is stable under the reflection in this third step.
\end{itemize}
Each of the three operations has an inverse, and thus  $\# \sapell_n \geq \tfrac{1}{8} \# \sap_n$. \qed

\subsection{Regulation global join polygons: three key properties}\label{s.threekey}

\begin{definition}\label{d.rgjp}
Let $k,\ell \in 2\N$
satisfy $k/2 \leq \ell \leq 35 k$. 
Let $\dubbub_{k,\ell}$ denote the set of {\em regulation global join} polygons (with length pair $(k,\ell)$), whose elements are formed by
Madras joining the polygon pair $\big(\phi^1,\phi^2 + \vec{u}\big)$, where
\begin{enumerate}
\item  $\phi^1 \in \sapell_k$;
\item $\phi^2 \in \sapr_\ell$;
\item $\ymax\big( \phi^2 + \vec{u} \big) \in \big[ y\big(\eastsouth(\phi^1)\big),y\big(\eastsouth(\phi^1)\big) +  \lfloor k^{1/2}/{10} \rfloor - 1 \big]$;
\item and $\big(\phi^1,\phi^2 + \vec{u} \big)$ is Madras joinable.
\end{enumerate}
Let $\dubbub$ denote the union of the sets $\dubbub_{k,\ell}$ over all such choices of $(k,\ell) \in 2\N \times 2\N$.
\end{definition}

Note that  Lemma~\ref{l.strongjoin} implies that whenever a polygon pair is joined to form an element of $\dubbub$, this pair is globally Madras joinable. 

The three key properties of regulation polygons are now stated as propositions.

\subsubsection{An exact formula for the size of $\dubbub_{k,\ell}$} 

\begin{proposition}\label{p.rgjpnumber}
Let $k,\ell \in 2\N$ satisfy $k/2 \leq \ell \leq 35 k$.  For any $\phi \in \dubbub_{k,\ell}$, there is a unique choice of 
$\phi^1 \in \sapell_k$, $\phi^2 \in \sapr_\ell$ and $\vec{u} \in \Z^2$ for which $\phi = J\big(\phi^1,\phi^2 + \vec{u} \big)$. 
We also have that
$$
\big\vert \dubbub_{k,\ell} \big\vert = \lfloor  k^{1/2}/{10} \rfloor \cdot  \big\vert \sapell_k \big\vert \, \big\vert \sapr_\ell \big\vert \, . 
$$
\end{proposition}

\subsubsection{The law of the initial path of a given length in a random regulation polygon}

We show that initial subpaths in regulation polygons are distributed independently of the length of the right polygon.

\begin{definition}
Write $\psapleft{n}$ for the uniform law on $\sapell_n$.
\end{definition}

\begin{proposition}\label{p.leftlong}
Let $k,\ell \in 2\N$ satisfy $k/2 \leq \ell \leq 35 k$. 
Let $j \in \N$ satisfy $j \leq k/2 - 1$. For any $\phi \in \sapell_k$,
$$
\psap_{k+\ell + 16} \Big( \Gamma_{[0,j]} = \phi_{[0,j]} \, \Big\vert \, \Gamma \in \dubbub_{k,\ell} \Big)
  =   \psapleft{k} \Big( \Gamma_{[0,j]} = \phi_{[0,j]} \Big) \, .
$$
\end{proposition}

\subsubsection{Regulation polygon joining is almost injective}

We reviewed in Section~\ref{s.polyjoin} how a form of~(\ref{e.polyjoin}) is obtained in~\cite{ThetaBound} by establishing efficient plaquette identification. We need something similar in our present endeavour. 
This next result is a counterpart to 
\cite[Proposition~4.2]{ThetaBound}.

\begin{proposition}\label{p.polyjoin.aboveonehalf}
For any $\Theta > 0$, there exist  $c > 0$ and $n_0 \in \N$ such that the following holds. Let  $i \geq 4$ be an integer and let $n \in 2\N \cap \big[ 2^{i+2},2^{i+3} \big]$ satisfy $n \geq n_0$.
Suppose that 
 $\mathsf{R}$
is a subset of $2\N \cap [2^i,2^{i+1}]$
that contains an element  $\kmac$ such that
$\max \big\{ \theta_{\kmac} , \theta_{n - 16 - \kmac} \big\} \leq \Theta$.
Then
$$
 \bigg\vert \, \bigcup_{j \in \mathsf{R}} \dubbub_{j,n-16-j} \, \bigg\vert \geq c \, \frac{n^{1/2}}{\log n} \sum_{j  \in \mathsf{R}} p_j \, p_{n-16-j}  \, .
$$
\end{proposition}

Section~\ref{s.pjprep} concludes with three subsections consecutively devoted to the proofs of these propositions.

\subsection{Proof of Proposition~\ref{p.rgjpnumber}} 
We need only prove the first assertion, the latter following directly. To do so, it is enough to determine the junction plaquette associated to the Madras join that forms $\phi$. 
By Lemma~\ref{l.strongjoin}, any such polygon pair $(\phi^1,\phi^2 + \vec{u})$ is globally Madras joinable. Recalling Definition~\ref{d.globallymj}, the associated junction plaquette is thus a global join plaquette of~$\phi$. 

By planarity, the complement of the union of the edges that comprise~$\phi$ is a bounded region in $\R^2$.
Also by planarity, the removal from the closure of this region of the closed unit square associated to any of $\phi$'s global join plaquettes will disconnect $\northeast(\phi)$ and $\eastsouth(\phi)$.
From this, we see that
\begin{itemize}
\item  each global join plaquette of $\phi$ has one horizontal edge traversed in the outward journey along $\phi$ from $\northeast(\phi)$ to $\eastsouth(\phi)$, and one traversed on the return journey;
\item and the set of $\phi$'s global join plaquettes is totally ordered under a relation in which the upper element is both reached earlier on the outward journey and later on the return.
\end{itemize}



We thus infer that the map that sends a global join plaquette $P$ of $\phi$ to the length of the polygon in $\phi \, \Delta \, P$ that contains $\northeast(\phi) = 0$ is injective. Since this length must equal $k+8$ for any admissible choice of 
$\phi^1 \in \sapell_k$, $\phi^2 \in \sapr_\ell$ and $\vec{u} \in \Z^2$, this choice is unique, and we are done. \qed


  
\medskip 

Note further that the set of join locations stipulated by the third condition in Definition~\ref{d.rgjp} is restricted to an interval of length of order~$k^{1/2}$. The restriction causes the formula in Proposition~\ref{p.rgjpnumber} to hold. It may be that many elements of~$\sapell_k$ and~$\sapr_\ell$ have heights much exceeding $k^{1/2}$, so that, for pairs of such polygons, there are many more than an order of  $k^{1/2}$ choices of translate for the second element that result in a globally Madras joinable polygon pair. The term {\em regulation} has been attached to indicate that a specific rule has been used to produce elements of $\dubbub$ and to emphasise that such polygons do not exhaust the set of polygons that we may conceive as being globally joined.

\subsection{Deriving Proposition~\ref{p.leftlong}}

The next result is needed for this derivation.

\begin{lemma}\label{l.leftlongbasic}
Consider a globally Madras joinable polygon pair $(\phi^1,\phi^2)$.
Let $j \in \N$ be such that $\phi^1_j$ equals $\southeast(\phi^1)$.  The join polygon $J(\phi^1,\phi^2)$ has the property that the initial subpaths $\phi^1_{[0,j-1]}$ and $J(\phi^1,\phi^2)_{[0,j-1]}$ coincide.
\end{lemma}
\noindent{\bf Proof.} The detailed structure of Madras joining from~\cite[Section 4.1]{ThetaBound} will be needed.
Since  $(\phi^1,\phi^2)$ is globally Madras joinable, the northeast vertex of the join polygon coincides with $\northeast(\phi^1)$. As such, it suffices to confirm that, absent its final edge, $\phi^1$'s left path $\phi^1_{[0,j-1]}$ forms part of $J(\phi^1,\phi^2)$.

Suppose that a vertex $v$ in $\phi^1$
is such that  the horizontal line segment extending to the right from~$v$ intersects the vertex set of~$\phi^1$ only at~$v$. Any such $v$
must equal~$\phi^1_i$ for some $i \geq j$ (i.e., must lie in the right path of $\phi^1$), because any  vertex of the left path of $\phi^1$, except its endpoint $\southeast(\phi^1)$, is directly to the left of some other vertex of~$\phi^1$.

 Review \cite[Figure 1]{ThetaBound}.
Applied to the present context, each sketch depicts a locale of $\phi^1$ in which surgery (and joining) will occur. A special vertex $Y$ is marked by a circle or a disc in each sketch. Recall from the Madras join construction that the horizontal line segment emanating rightwards from $Y$ has no intersection with the vertex set of $\phi^1$, except possibly at its leftmost point;
moreover, this is also true when $Y$ is replaced by either $Y - e_2$ or $Y + e_2$.

\begin{figure}
    \begin{center}
      \includegraphics[width=0.7\textwidth]{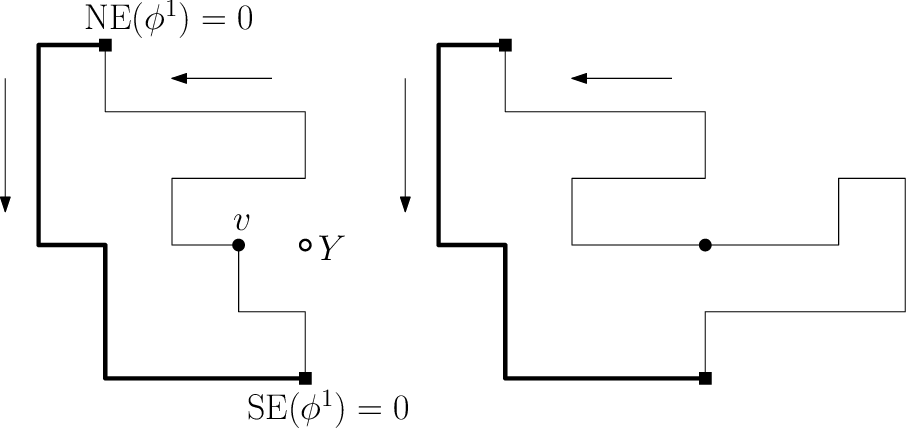}
    \end{center}
    \caption{A case in the proof of Lemma~\ref{l.leftlongbasic}. The polygon $\phi^1$ depicted on the left lies in $\sap_n$ for some $n \in 2\N$. The polygon is indexed counterclockwise from $\northeast(\phi^1)$ as the arrows indicate.
   Recalling \cite[Figure 1]{ThetaBound},  
   the circle in the left sketch, and the disk in the right, marks the point $Y$ used in the formation of the right sketch $\modify{\tau}$ with $\tau = \phi^1$, (during the process of Madras joining with an undepicted~$\phi^2$). 
    Case~IIcii obtains for $\phi^1$. The left sketch disk marks a vertex $v$ as discussed in the proof.  Emboldened is the left path of $\phi^1$, which passes undisturbed by the surgery from the left sketch to the right to form an initial subpath of the postsurgical polygon.}\label{f.caseiicii}
  \end{figure}

Thus, we see that, whichever case for~\cite[Figure 1]{ThetaBound} obtains, any vertex of $\phi^1$ lying in $\big\{ Y-e_2,Y,Y + e_2 \big\}$ has the form $\phi^1_i$ for $i \geq j$. Inspection of the figure shows that the surgery conducted to form the Madras join polygon leaves $\phi^1_{[0,j-1]}$ unaltered: in case IIb, this inference depends on further noting that both $Y- e_1$ and $Y - e_1 + e_2$ take the form $\phi^1_i$ for $i \geq j$; in case IIci, it depends on noting that $Y - e_2 - e_1$ and $Y - e_1$ take such a form. 
(Case IIcii is illustrated in Figure~\ref{f.caseiicii}.) Moreover, this assertion of constancy under surgery  includes cases IIIa, IIIb, IIIci and IIIcii which are not depicted but are discussed in the caption of~\cite[Figure 1]{ThetaBound}. 

In fact, in all cases other than case IIa,
even $\phi^1_{[0,j]}$ suffers no alteration in surgery.
It is case IIa and the possibility that $Y + e_2 = \southeast(\phi_1)$ that leads us to claim merely that 
$\phi^1_{[0,j-1]} = J(\phi^1,\phi^2)_{[0,j-1]}$ in 
Lemma~\ref{l.leftlongbasic}. \qed

\medskip

\noindent{\bf Proof of Proposition~\ref{p.leftlong}.} By the first assertion in Proposition~\ref{p.rgjpnumber}, we find that the conditional distribution of 
$\psap_{k+\ell+16}$ given that~$\Gamma \in \dubbub_{k,\ell}$ has the law of the output in this procedure:
\begin{itemize}
\item select a polygon according to the law $\psapleft{k}$;
\item independently select another uniformly among elements of $\sapr_\ell$;
\item choose one of the $\lfloor \litnum \, k^{1/2} \rfloor$ vertical displacements specified in Definition~\ref{d.rgjp} uniformly at random, and shift the latter polygon by this displacement;
\item translating the displaced polygon by a suitable horizontal shift, Madras join the first and the translated second polygon to obtain the output.
\end{itemize} 
Denote by $\phi^1$ the element of $\sapell_k$ selected in the first step. This polygon is left-long; let $m \geq k/2$ satisfy $\phi^1_m = \southeast(\phi^1)$. By Lemma~\ref{l.leftlongbasic}, the output polygon has an initial subpath that coincides with the left path $\phi^1_{[0,m-1]}$ of~$\phi^1$ (absent the final edge). Since $j \leq k/2 - 1 \leq m - 1$, we see that the initial length~$j$ subpath of the output polygon coincides with $\phi^1_{[0,j]}$. Thus, its law is that of $\Gamma_{[0,j]}$ where the polygon $\Gamma$ is distributed according to~$\psapleft{k}$. \qed

\subsection{Proof of Proposition~\ref{p.polyjoin.aboveonehalf}}

This argument is based on the proof of \cite[Lemma~4.12]{ThetaBound}, (which is the principal component of Proposition~4.2 in that paper).
We mentioned in Section~\ref{s.polyjoin} that efficient plaquette identification is proved to work in~\cite{ThetaBound}
only when the length of the joined polygon lies in the index set  $\highpolynum_\maceta$ (for a suitable choice of $\maceta > 0$).
It is Corollary~\ref{c.globaljoin} that demonstrates that identification is efficient. We use this approach again now.

Let $\Theta > 0$, $i \geq 4$ and $n \in 2\N \cap [2^{i+2},2^{i+3}]$. Suppose that a set $\mathsf{R}$ contains an element $k^*$ in the way that Proposition~\ref{p.polyjoin.aboveonehalf} describes.
 We begin by pointing out that 
  for any choice of $\maceta$ exceeding $2\Theta$,
  the assertion that $n \in \highpolynum_\maceta$ holds provided that $n$ is at least the value $n_0(\maceta) = \mu^{16/(\maceta - 2\Theta)}$. 
  That this is true follows from 
the straightforward $p_n \geq p_{n-16}$, polygon superadditivity $p_{n-16} \geq p_{\kmac} p_{n-16-\kmac}$ (see~\cite[(3.2.3)]{MS93}) and the hypothesised $\max \big\{ \thet_{\kmac},\thet_{n-16-\kmac} \big\} \leq \Theta$. 

We now consider  a multi-valued map $\Psi:A \to \mathcal{P}(B)$. 
We take
$A = \bigcup_{j \in \mathsf{R}} \, \sapell_j \times \sapr_{n-16-j}$  and   $B =  \cup_{j \in \mathsf{R}} \,  \dubbub_{j,n-16-j}$ (which is a subset of $\sap_n \cap \dubbub$).
Note  the indices concerned in specifying $A$ are non-negative  because $\max \mathsf{R} \leq 2^{i+1}$, $n \geq 2^{i+2}$ and $i \geq 3$ ensure that $n - 16 -j \geq 0$ whenever $j \in \mathsf{R}$. (Moreover,  the choice $(k,\ell) = (j,n-16-j)$ in Definition~\ref{d.rgjp} does satisfy the bounds stated there. We use $i \geq 4$ to verify that $k/2 \leq \ell$.)
We specify $\Psi$ so that  the generic domain point 
 $(\phi^1,\phi^2) \in \sapell_j \times \sapr_{n-16-j}$, $j \in \mathsf{R}$, is mapped to  
 the collection of length-$n$ polygons 
formed by Madras joining $\phi^1$ and $\phi^2 + \vec{u}$ as $\vec{u}$ ranges over the subset of $\strongjoin_{(\phi^1,\phi^2)}$ specified for the given pair $(\phi^1,\phi^2)$ by conditions $(3)$ and $(4)$ in Definition~\ref{d.rgjp}. (That this set is indeed a subset of $\strongjoin_{(\phi^1,\phi^2)}$ is noted after Definition~\ref{d.rgjp}.)

We have that $\big\vert  \Psi\big( (\phi^1,\phi^2) \big) \big\vert = \lfloor j^{1/2}/10 \rfloor$. 
 Since $j \geq 2^{i}$ and $n \leq 2^{i+3}$, we have that $j \geq  2^{-3}n$, so that $\big\vert  \Psi\big( (\phi^1,\phi^2) \big) \big\vert \geq  \lfloor \tfrac{1}{10} 2^{-3/2} n^{1/2}  \rfloor$.

Applying Lemma~\ref{l.polysetbound}, we learn that the number of arrows in $\Psi$ is at least
\begin{equation}\label{e.totalarrow.late}
 2^{-4} \cdot  \lfloor \tfrac{1}{10} 2^{-3/2} n^{1/2}  \rfloor \sum_{j \in \mathsf{R}} p_j p_{n-16-j}  \, .
\end{equation}

For a constant $\Ch > 0$, denote by 
$$
H_n = H_n(\Ch) = \Big\{ \phi \in B : \big\vert \globaljoin_\phi \big\vert \geq  \Ch \log n \Big\} 
$$
the set of length-$n$ polygons with a {\em high} number of global join plaquettes.

We now set the value of $\Ch > 0$.
Specify a parameter $\maceta = 2\Theta + 1$, and recall that 
$n \in \highpolynum_\maceta$ is known provided that $n \geq n_0(\maceta)$. We set the parameter $\Ch$ equal to the value  $\Cgjtwo$ furnished by Corollary~\ref{c.globaljoin} applied with $\Cgjone$ taking the value $2\Theta + 1$ and for the given value of 
$\maceta$ (which happens to be the same value). This use of the corollary also provides a value of $\Cgjthree$, and,  
since $\psap_n \big( \vert \globaljoin_\Gamma \vert \geq \Ch \log n \big) = p_n^{-1} \vert H_n \vert$, we learn that
$\vert H_n \vert \leq  \Cgjthree n^{-(2\Theta + 1)}  p_n$.

There are two cases to analyse, the first of which is specified by at least one-half of the arrows in $\Psi$ pointing to $H_n$.
In the first case, we have
\begin{eqnarray*}
 & &  n \times  \Cgjthree n^{-(2\Theta + 1)}  p_n  \, \geq \,  \max \big\{ \vert \Psi^{-1}(\phi) \vert: \phi \in H_{n} \big\} \times  \big\vert H_{n} \big\vert 
  \, \geq \, \tfrac{1}{2} \# \big\{ \,  \textrm{arrows in $\Psi$} \, \big\} \\
  & \geq & 
 2^{-5} \lfloor \tfrac{1}{10} 2^{-3/2} n^{1/2} \rfloor \sum_{j \in \mathsf{R}} p_j p_{n-16-j} \, \geq \, \alpha \, n^{1/2} p_{\kmac} p_{n-16- \kmac} \, \geq \, \alpha \, n^{1/2} n^{-2\Theta} \mu^{n -16} \, .  
\end{eqnarray*} 
where the value of $\alpha$ may be chosen to equal $\tfrac{1}{10} 2^{-7}$ since $n$ is supposed to be least $n_0$, and the latter quantity may be increased if need be.
The consecutive inequalities are due to: $\big\vert \Psi^{-1}(\phi) \big\vert \leq n$ for  $\phi \in \sap_n$; the first case obtaining; the lower bound
(\ref{e.totalarrow.late}); the membership of $\mathsf{R}$ by the hypothesised index $\kmac$; and the hypothesised properties of $\kmac$. We have found that $p_n \geq c \, n^{1/2} \mu^n$ with $c = \mu^{-16} \alpha \Cgjthree^{-1}$. Since $p_n \leq \mu^n$ (as noted in the introduction), the first case ends in contradiction if 
we stipulate that the hypothesised lower bound $n_0$ on $n$ is at least $c^{-2} = 100 \cdot 2^{14} \mu^{32} \Cgjthree$ (where this bound has a $\Theta$-dependence that is transmitted via $\Cgjthree$).

The set of preimages under $\Psi$ of a given 
$\phi \in B \subset \sap_n$ may be indexed by the junction plaquette 
associated to the Madras join polygon $J(\phi^1,\phi^2)$ that equals $\phi$. Recalling Definition~\ref{d.globallymj}, this plaquette is a global join plaquette of $\phi$. Thus,  
\begin{equation}\label{e.phisap}
\vert \Psi^{-1}(\phi) \vert \leq \vert \globaljoin_\phi \vert
 \, \, \,
\textrm{for any $\phi \in B$} \,  
.
\end{equation}

We find then that, when the second case obtains,
$$
 \vert B \vert \geq  \big\vert B \setminus H_n \big\vert  \geq    \Ch^{-1} \big( \log n \big)^{-1} \times
  2^{-5} \lfloor \tfrac{1}{10} 2^{-3/2} n^{1/2} \rfloor  \sum_{j \in \mathsf{R}} p_j p_{n-16-j}  \, .
$$
Recalling the definition of the set $B$ completes the proof of Proposition~\ref{p.polyjoin.aboveonehalf}. \qed

\section{The snake method: general elements}\label{s.five}

In this section, we recall from~\cite{snakemethodpattern} the snake method. 
The method is used to prove upper bounds on the closing probability, and assumes to the contrary that to some degree this probability has slow decay. For the technique to be used, two ingredients are needed.
\begin{enumerate}
\item A charming snake is a walk or polygon $\ga$ many of whose subpaths beginning at $\northeast(\ga)$  have high conditional closing probability, when extended by some common length. It must be shown that charming snakes are not too atypical.
\item A general procedure is then specified in which 
a
charming snake is used to manufacture huge numbers of alternative self-avoiding walks. These alternatives overwhelm the polygons in number and show that the closing probability is very small, contradicting the assumption.
\end{enumerate}

The first step certainly makes use of the assumption of slow decay on the closing probability. It is not however a simple consequence of this assumption. After this section, we will carry out the first step via the polygon joining apparatus we have constructed to prove Theorem~\ref{t.thetexist}; in~\cite{snakemethodpattern},
the first step was implemented via a technique of Gaussian pattern fluctuation to prove 
 the general dimensional Theorem~\ref{t.closingprob}.

This section is devoted to explaining the second step, which is formulated as a general tool, valid in any dimension $d \geq 2$.  

\subsection{Some necessary notation and tools}

To set up the apparatus of the method, we present some preliminaries.

\subsubsection{The two-part decomposition}

In the snake method, we represent any given walk $\gamma$
in a {\em two-part} decomposition. This consists of an ordered pair of walks $(\gamma^1,\gamma^2)$ that emanate from
a certain common vertex and that are disjoint except at that vertex.
The two walks are called the {\em first part} and the {\em second part}.
To define the decomposition, consider any walk $\gamma$ of length~$n$. We first mention that the common vertex is chosen to be the northeast vertex $\northeast(\gamma)$.
Choosing $j \in [0,n]$  
so that $\gamma_j = \northeast(\gamma)$, the walk $\gamma$ begins at $\ga_0$ and approaches $\northeast(\ga)$ along the subwalk $\gamma_{[0,j]}$,
and then continues to its endpoint $\ga_n$ along the subwalk $\gamma_{[j,n]}$. The reversal $\reverse\gamma_{[0,j]}$ of the first walk, and the second walk $\gamma_{[j,n]}$, form a pair of walks that emanate from $\northeast(\gamma)$. (When $j$ equals zero or $n$,
one of the walks is the length zero walk that only visits $\northeast(\gamma)$; in the other cases, each walk has positive length.) The two walks will be the two elements in the two-part decomposition; all that remains is to order them, to decide which is the first part. 
If one of the walks visits $\northeast(\gamma) - e_1$ at its first step, then it is set equal to $\gamma^1$; of course, $\gamma^2$ is then the other walk. If neither visits  $\northeast(\gamma) - e_1$, then one of the walks has length zero. This walk is $\gamma^1$ in this case.

We use square brackets to indicate the two-part decomposition, writing $\ga = [\ga^1,\ga^2]$.


As a small aid to visualization, it is useful to note that if the first $\ga^1$ part of a two-dimensional walk $\ga$ for which $\northeast(\gamma) = 0$ has length~$j$, then $\gamma^1:[0,j] \to \Z^2$ satisfies 
\begin{itemize}
\item $\gamma^1_0 = 0$ and $\gamma^1_1 = - e_1$;
\item $y \big( \gamma^1_i \big) \leq 0$ for all $i \in [0,j]$;
\item $\gamma^1_i \not\in \N \times \{ 0 \}$ for any $i \in [1,j]$.
\end{itemize}

  \begin{figure}
    \begin{center}
      \includegraphics[width=0.3\textwidth]{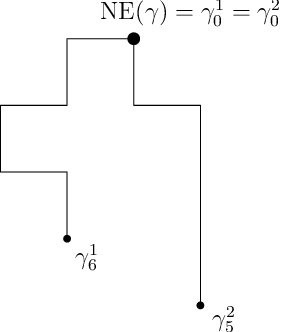}
    \end{center}
    \caption{The two-part decomposition of a walk of length eleven.
}\label{f.firstpart}
  \end{figure}



\subsubsection{Notation for walks not beginning at the origin}\label{s.sawfree}
Let $n \in \N$. We write $\sawfree_n$ for the set of self-avoiding walks~$\ga$ of length~$n$ (without stipulating the location~$\ga_0$). We further write $\saw_n^0$ for the subset of $\sawfree_n$ whose elements~$\ga$ verify $\northeast(\ga) = 0$. Naturally, an element $\ga \in \saw_n^0$ is said to close (and be closing) if $\dist \ga_n - \ga_0 \dist = 1$. The uniform law on $\saw_n^0$ will be denoted by $\psaw_n^0$. The sets $\saw_n^0$ and $\saw_n$ are in bijection via a clear translation; we will use this bijection implicitly.

\subsubsection{First parts and closing probabilities}

First part lengths with low closing probability are rare, as the following result, \cite[Lemma 2.2]{snakemethodpattern}, attests.

\begin{lemma}\label{l.closecard}
Let $n \in 2\N + 1$ be such that, for some $\alphamac > 0$, $\psaw_n \big( \Gamma \closes \big) \geq n^{-\alphamac}$. For any $\deltamac > 0$,  the set of $i \in [0,n]$ for which
$$
 \# \Big\{ \gamma \in \saw^0_n: \vert \gamma^1 \vert = i \Big\}
 \geq n^{\alphamac + \deltamac}  \cdot \# \Big\{ \gamma \in \saw^0_n: \vert \gamma^1 \vert = i  \, , \, \gamma \closes \Big\}
$$
has cardinality at most $2 n^{1- \deltamac}$. 
\end{lemma}

\subsubsection{Possible first parts and their conditional closing probabilities}\label{s.possparts}

For $n \in \N$, let $\fpart_n \subseteq \saw_n$ denote  the set of walks  $\gamma:[0,n] \to \Z^2$
whose northeast vertex is $\gamma_0 = 0$. 
We wish to view $\fpart_n$ as the set of possible first parts of walks $\phi \in \saw^0_m$ of some length $m$ that is at least $n$.
(We could in fact be more restrictive in specifying $\fpart_n$, stipulating if we wish that any element $\ga$ satisfies $\ga_1 = -e_1$. As it stands, $\fpart_n$ actually contains all possible second parts. What matters, however, is only that $\fpart_n$ contains all possible first parts.)

\begin{figure}
    \begin{center}
      \includegraphics[width=0.4\textwidth]{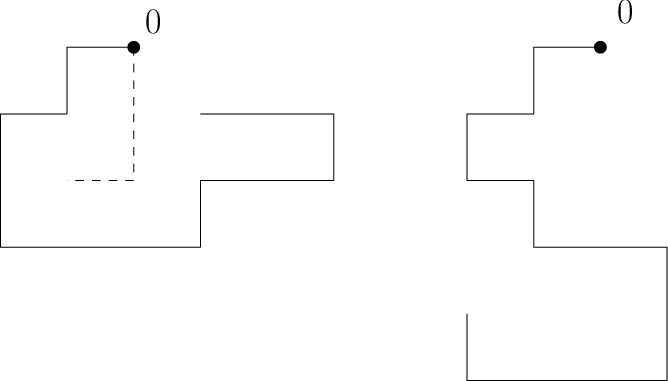}
    \end{center}
    \caption{{\em Left}: the bold $\phi \in \saw_{14}$ and dashed $\ga \in \saw_3$ are such that $[\phi,\ga]$ is a two-part decomposition. Note that $\phi \in \fpart_{14,14+3} \cap \fpart_{14,14+4}^c$. {\em Right}: An element of $\cap_{m = 1}^\infty\fpart_{14,14+m}$.}\label{f.twophi}
  \end{figure}

Note that, as Figure~\ref{f.twophi} illustrates, for given $m > n$, only some elements of $\fpart_n$ appear as such first parts, and we now record notation for the set of such elements. Write $\fpart_{n,m} \subseteq \fpart_n$ for the set of $\gamma \in \fpart_n$ for which there exists an element $\phi \in \saw_{m-n}$ (necessarily with $\northeast(\phi) = 0$) such that $[\gamma,\phi]$ is the two-part decomposition of some element $\chi \in \saw_m$ with $\northeast(\chi) = 0$. 

 In this light, we now define the conditional closing probability 
$$
 q_{n,m}:\fpart_{n,m} \to [0,1] \,   , \, \, \,  \, q_{n,m}(\gamma) = \psaw^0_m \Big( \Gamma \closes \, \Big\vert \, \Gamma^1 = \gamma \Big) \, ,
$$ 
where here $m,n \in \N$ satisfy $m > n$; note that since $\gamma \in \fpart_{n,m}$, the event in the conditioning on the right-hand side occurs for some elements of $\saw_m$, so that the right-hand side is well-defined. 

We also identity a set of first parts with {\em high} conditional closing probability: for $\alpha > 0$, we write
$$
 \hPhi_{n,m}^\alpha = \Big\{ \gamma \in \fpart_{n,m}: q_{n,m}(\gamma) > m^{-\alpha} \Big\} \, .
$$

\subsection{The snake method: recall of its general apparatus}

\subsubsection{Parameters}

The snake method has three exponent parameters:
\begin{itemize}
\item the inverse charm $\alpha > 0$;
\item the snake length $\beta \in (0,1]$;
\item and the charm deficit $\esm \in (0,\beta)$.
\end{itemize}
It has two index parameters:
\begin{itemize}
\item $n \in 2\N + 1$ and $\ell \in \N$, with $\ell \leq n$.
\end{itemize}

\subsubsection{Charming snakes}

Here we define these creatures.

\begin{definition}
Let $\alpha > 0$, $n \in 2\N + 1$,  $\ell \in [0,n]$,
$\gamma \in \fpart_{\ell,n}$, and $k \in [0,\ell]$ with $\ell - k \in 2\N$.  We say that $\ga$ is $(\alpha,n,\ell)$-{\em charming} at (or for) index $k$ if 
\begin{equation}\label{e.charming}
\psaw^0_{k + n - \ell} \Big( \Gamma \closes \, \Big\vert \, \vert \Gamma^1 \vert = k \, , \, \Gamma^1 = \gamma_{[0,k]} \Big) > n^{-\alpha} \, .
\end{equation}
\end{definition}
(The event  that $\vert \Gamma^1 \vert = k$ in the conditioning is redundant and is recorded for emphasis.)
Note that an element $\ga \in \fpart_{\ell,n}$  is $(\alpha,n,\ell)$-charming at  index~$k$ if in selecting a length~$n-\ell$ walk beginning at $0$ uniformly among those choices that form the second part of a walk whose first part is $\ga_{[0,k]}$, the resulting $(k + n - \ell)$-length walk closes with probability exceeding $n^{-\alpha}$. 
(Since we insist that $n$ is odd and that $\ell$ and $k$ share their parity, the length $k + n - \ell$ is odd; the condition displayed above could not possibly be satisfied if this were not the case.)
Note that, for any $\ell \in [0,n]$,  $\ga \in \fpart_{\ell,n}$ is $(\alpha,n,\ell)$-charming at the special choice of index $k = \ell$ precisely when $\ga \in \hPhi_{\ell,n}^{\alpha}$. When $k < \ell$ with $k+n - \ell$ of order~$n$, the condition that 
 $\ga \in \fpart_{\ell,n}$ is $(\alpha,n,\ell)$-charming at  index $k$ is {\em almost} the same as  $\ga_{[0,k]} \in \hPhi_{k,k + n - \ell}^{\alpha}$; (the latter condition would involve replacing $n$ by $k + n -\ell$ in the right-hand side of~(\ref{e.charming})). 

For $n \in 2\N + 1$, $\ell \in [0,n]$, $\alpha,\beta > 0$ and $\esm \in (0,\beta)$, define the charming snake set 
\begin{eqnarray*}
\mathsf{CS}_{\beta,\esm}^{\alpha,\ell,n}  & =  & \Big\{ \ga \in \fpart_{\ell,n}: \textrm{$\ga$ is $(\alpha,n,\ell)$-charming} \\
 & & \quad \textrm{for at least $n^{\beta- \esm}/4$ indices belonging to the interval $\big[ \ell - n^{\beta}, \ell \big]$} \Big\}
 \, .
\end{eqnarray*}
A charming snake is an element of 
$\mathsf{CS}_{\beta,\esm}^{\alpha,\ell,n}$ and, as such, it is a walk. To a charming snake $\gamma$, however, 
is associated the sequence 
$\big( \gamma_{[0,\ell - n^{\beta}]},\ga_{[0,\ell - n^{\beta} + 1]}, \cdots , \ga_{[0,\ell]} \big)$ of $n^{\beta} + 1$ terms. 
This sequence may be depicted as evolving as an extending snake does; since $\gamma \in \mathsf{CS}_{\beta,\esm}^{\alpha,\ell,n}$, the sequence has many $(\alpha,n,\ell)$-charming elements, 
for each of which there is a high conditional closing probability for extensions of a {\em shared} continuation length~$n - \ell$.

\subsubsection{A general tool for the method's second step}

For the snake method to produce results, we must work with a choice of parameters for which $\beta - \esm - \alpha > 0$. (The method could be reworked to handle equality in some cases.) Here we present the general tool for carrying out the method's second step. 
The tool asserts that, if $\beta - \esm - \alpha > 0$ and even a very modest proportion of snakes are charming, then the closing probability drops precipitously.

\begin{theorem}\label{t.snakesecondstep}
Let $d \geq 2$. Set $c = 2^{\tfrac{1}{5(4d + 1)}} > 1$ 
and set $K = 20 (4d+1) \tfrac{\log(4d)}{\log 2}$.
Suppose that the exponent parameters satisfy $\delta = \beta - \esm - \alpha > 0$.
If the index parameter pair $(n,\ell)$ satisfies $n \geq K^{1/\delta}$ and 
\begin{equation}\label{e.afewcs}
 \psap_{n+1} \Big( \Gamma_{[0,\ell]} \in \mathsf{CS}_{\beta,\esm}^{\alpha,\ell,n}  \Big)  \geq c^{- n^{\delta}/2} \, ,
\end{equation}
then
$$
\psaw_n \Big( \Gamma \closes \Big) \leq 2 (n+1) \,  c^{-n^{\delta}/2} \, .
$$ 
\end{theorem}
The proof appears in~\cite[Section 3.2]{snakemethodpattern}.

Note that since the closing probability is predicted to have polynomial decay, the hypothesis~(\ref{e.afewcs}) is never satisfied in practice. For this reason, the snake method will always involve argument by contradiction, with~(\ref{e.afewcs}) being established under a hypothesis that the closing probability is, to some degree, decaying slowly.

\section{An impression of the ideas for the main proof}\label{s.expo}

We have assembled the principal inputs needed to prove Theorem~\ref{t.thetexist} and will begin proving this result in the next section. The proof is rather technical and delicate, and, before we begin it, we attempt to set the reader's bearings by giving a necessarily impressionistic account of the ideas that drive the derivation.

Remember that our basic aim in this paper to prove that, when $d=2$, the closing probability has a rate of decay more rapid than $n^{-1/2 - \eps}$, for some $\eps > 0$.
Since we will proceed by proof by contradiction, it is helpful to reflect on the circumstance that the opposite case holds, which we may summarise by the informal assertion that $\psaw_n (\Ga \closes) \geq n^{-1/2 - o(1)}$. In this section, we will refer to this bound, assumed for the present purpose of exposition to hold for all high odd $n$, as the {\em standing assumption}. Recalling~(\ref{e.closepc})
and $c_n \geq \mu^n$,
our assumption forces $p_n \geq n^{-3/2 - o(1)} \mu^n$ for even $n$. (The opposing bound $p_n \leq n^{-3/2 + o(1)} \mu^n$ holds for typical such $n$ by Theorem~\ref{t.polydev}. The standing assumption thus yields the perhaps provocatively implausible bound $c_n \leq n^{o(1)} \mu^n$ for typical $n$. This remark is however an aside, because we will make no use of this bound.)

The lower bound on polygon number forces all sufficiently high even numbers to be high polygon number indices in the sense of Definition~\ref{d.highpolynum}, so that the hypothesis $n \in \highpolynum_\zeta$ of Corollary~\ref{c.globaljoin} (for any given $\zeta > 3/2$) is  verified for high even $n$. This means that the circumstance of efficient plaquette identification obtains, so that a form of~(\ref{e.polyjoin}) holds for large even~$n$.

Now, the right-hand side of~(\ref{e.polyjoin}) is an expression for the number of ways of joining together  pairs of polygons each of size of order $n$, provided that a set of order~$n^{1/2}$ possibilities is used for each pair of polygons. More joinings may be possible, but we do not know this for sure: we will call the Madras join polygons arising from a specific set of order $n^{1/2}$ join locations `regulation global join polygons' in the upcoming proof, though here we prefer to call them simply `global join polygons', to emphasise their key geometric feature of being comprised of two chambers of the same order of length.

Now a very distinctive geometric feature obtains under the standing assumption: a non-negligible fraction (of size at least $n^{-o(1)}$) of polygons of length a given high even value $n$ must be global join polygons.
We may call this feature `ample supply of global join polygons'. The reason that such polygons are in ample supply is that, if these special polygons represented only a negligible fraction of polygons of a given high length, then~(\ref{e.polyjoin}) could be improved by the presence of an addition factor of $n^\eps$
on its right-hand side, and the inference that $p_n \mu^{-n} \leq n^{-3/2 - \eps}$ would result, contrary to the bound $p_n \geq n^{-3/2 - o(1)} \mu^n$ that has been seen to follow from the standing assumption.

The ample supply of global join polygons will play a critical role in validating the snake method hypothesis~(\ref{e.afewcs}). 
Under the standing assumption, we will endeavour to verify this hypothesis with the parameter choice $\beta =1$; this choice reflects the ambition to implement the method in its strongest form, in which a charming snake whose length is macroscopic (of order $n$) must be shown to be not atypical. We further set, for reasons that will shortly be discussed, $\alpha = 1/2 + o(1)$ and $\eta = o(1)$.
(The bound that $\delta = \beta - \eta - \alpha$ be positive is satisfied; when we turn later to rigorous argument, the consideration that $\delta > 0$ dictates the degree to which we succeed in showing that the closing probability decays more quickly than $n^{-1/2}$.)

 To explain the basic meaning of the method hypothesis~(\ref{e.afewcs}),
 first recall that we have been treating $n$ in the preceding paragraphs as a general parameter and thus permitting it to be even. In the snake method, however, it is fixed at given high odd value, (reflecting its role as the length of a closing walk).
So fixing $n$, we also make a choice of the parameter $\ell \leq n$ of the same order as $n$. Think of a journey travelled counterclockwise along a typical polygon under $\psap_{n+1}$.
To each moment of the journey --  to each vertex of the polygon -- corresponds the first part of the walk formed by the removal from the polygon of the edge over which the traveller will next pass.  Since we take $\beta =1$, the journey includes order $n$ steps, so that the traveller passes over a macroscopic part of the polygon. To each first part, we consider a second part of fixed length $n-\ell$, chosen uniformly among the choices of second part that are admissible given the first part. The hypothesis~(\ref{e.afewcs})
is satisfied when it is with at least modest probability (for simplicity, say at least $n^{-o(1)}$), that a non-negligible fraction of these first parts (specifically, at least $n^{1 - o(1)}$ of them) are charming in the sense that this random selection of second part has a probability at least $n^{-1/2 - o(1)}$ of closing the first part, i.e., of ending next to the non-origin endpoint of the given first part. 

Roughly speaking, this choice of cutoff  $n^{-1/2 - o(1)}$, dictated by setting $\alpha = 1/2 + o(1)$,  is compatible with the standing assumption that the closing probability decays at rate very close to $n^{-1/2}$.
Indeed, if we were to vary the form of the hypothesis~(\ref{e.afewcs}) a little (calling the new version the `simple' form of the hypothesis), we could readily deduce that many first parts have the requisite property. In order to specify the simple form of hypothesis~(\ref{e.afewcs}),  we adopt the definition   that,   for a polygon of length $n+1$, a first part, of a length $k$ of order $n/2$, is `simply' charming if, when a second part extension of length $n -k$ is sampled, it closes the first part with probability at least $n^{-1/2 - o(1)}$.

Because the closing walks so formed correspond to polygons of length $n$, it is now a direct consequence of our standing assumption that, for a typical polygon sampled from $\psap_{n+1}$, most first parts are simply charming. Indeed, the number of first parts that fail to be simply charming is bounded above by Lemma~\ref{l.closecard}. (The proof of this result appears in~\cite{snakemethodpattern}, but it  is very straightforward.) That is, simply charming first parts are typical, with $n^{1 - o(1)}$ exceptions permitted among order $n$ lengths of first part.

Our discussion has led us to a clear expression of the challenge we face in implementing the snake method. The hypothesis~(\ref{e.afewcs}) is validated in its `simple' form, that is,  when we are prepared to decrease the second part length as the first part length rises (so that the sum of lengths equals $n$).
 In its actual form, of course,  we must accept that the second part length remains fixed (at the value $n-\ell$) even as the first part length $k$ increases
over a range of order $n$ values.  In fact, the simple version of the hypothesis has been comfortably validated. From the form of~(\ref{e.afewcs}), we see that
actually charming first parts must be shown to represent a non-negligible, $n^{-o(1)}$, fraction of all first parts for a non-negligible fraction of all length $n+1$ polygons; but first parts have been seen to be simply charming apart from a small, $n^{1-o(1)}$-sized, set of violations, with probability at least $1 - n^{-o(1)}$.

It is the ample supply of global join polygons that allows us to bridge the gap between the simple and actual forms of the snake method hypothesis. Building this bridge is a matter of understanding a certain similarity of measure between the polygon laws $\psap_{n+1}$ and $\psap_{m+1}$, where $m$ is odd and of order $n$. We are trying to gain an understanding that there is at least probability $n^{-o(1)}$ that a first part of a given length $k$ arising from the law~$\psap_{n+1}$ will have the property that a length $n - \ell$
second part will close with probability at least $n^{-1/2 - o(1)}$. By considering the choice $m = n - \ell + k$, we have this information for typical such~$k$ in the case that the polygon law is $\psap_{m+1}$ in place of $\psap_{n+1}$: this is the deduction made from knowing that the snake method hypothesis is valid in its simple form.

To build the bridge, recall that a non-negligible, $m^{-o(1)}$, fraction of length-$(m+1)$ polygons are global join polygons. There is an attractive random means of mapping such polygons into their length-$(n + 1)$ counterparts. Any such polygon is comprised of a left polygon and a right polygon. Pull out the right polygon, replace it by a polygon, uniformly chosen among the space of polygons whose length exceeds the removed polygon's by $n-m$ (whatever the sign of this difference), and rejoin the new polygon in one of the order $n^{1/2}$ spots that we regard as admissible for the formation of globally joined polygons. 

The resulting polygon has length $n+1$; (clearly, we are neglecting the details of the Madras join construction here). Moreover, this random function, mapping the space of length-$(m+1)$ globally joined polygons to the space of their length-$(n+1)$ counterparts, is roughly measure-preserving when projected to first parts that do not reach the location where surgery occurs (and thus are undisturbed by surgery).  

The probability that a first part of given length $k$ fails the simple charm test, when polygon length equals $m + 1 = n - \ell + k + 1$,
is low enough that it is much smaller than  the proportion of length-$(m+1)$ polygons that are global join polygons. (Ample supply of global join polygons means that the latter proportion, while not of unit order, is non-negligible. When we reprise these arguments rigorously, parameter choices may be made, in correspondance to different $o(1)$ terms, so that the failure probability for the simple charm test for a given length first part will indeed be much smaller than the global join polygon proportion.)

\begin{figure}
  \begin{center}
    \includegraphics[width=1.0\textwidth]{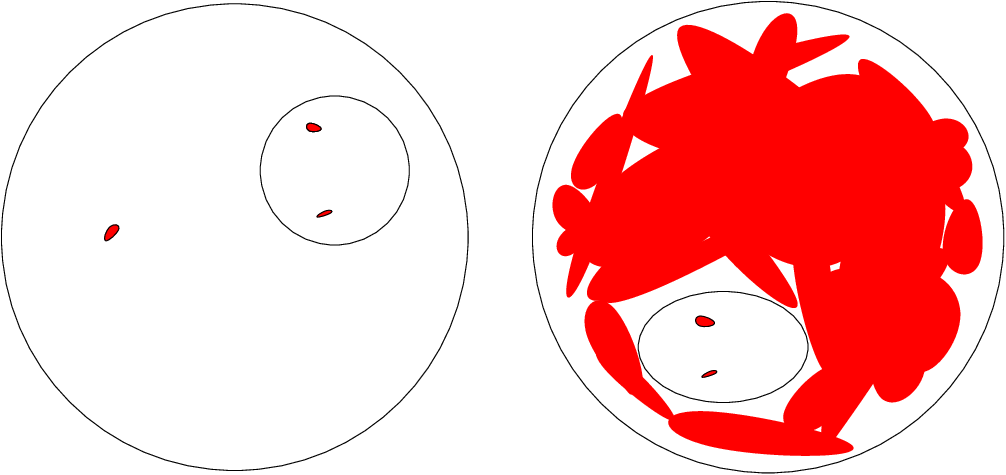}
  \end{center}
  \caption{For $k$ given of order $n$, the two large balls respectively represent the laws $\psap_{m+1}$ and $\psap_{n+1}$, where $m = k + n -\ell$. Imagine that each ball is partitioned into tiny disjoint cells that represent the polygons of the length in question. Inside the two balls, the cells associated to polygons whose length-$k$ first part fails to be charming for a second part extension of length $n - \ell$ are depicted in red. The smaller ball on the left represents  globally joined polygons of length~$m+1$. These are transported in a roughly area preserving way to the interior of the ellipse on the right, which represents globally joined polygons of length $n+1$. Under the standing assumption, the red region on the left has smaller area than does the ball, so that much of the ellipse interior on the right is seen to be not red. Moreover, the ellipse interior on the right represents a non-negligible area of the right-hand ball. Thus, even though the right-hand ball may be predominately red outside the ellipse, a non-negligible region on the right is not red, so that length~$k$ first parts are not atypically charming under~$\psap_{n+1}$.}\label{f.randomtransport}
\end{figure}

Consider now the random mapping of globally joined polygons of length $m+1$ into their length $n+1$ counterparts (and consult Figure~\ref{f.randomtransport} for an illustration of the argument). 
This map roughly preserves measure, so that the  charm evinced by a length-$k$ first part (namely, that a length-$(n-\ell)$ second part closes with probability at least $n^{-1/2 - o(1)}$)
is seen to be typical not only under a uniformly chosen globally joined polygon of length $m+1$
but also under the counterpart measure for length $n+1$. By the ample supply of globally joined polygons, a non-negligible fraction of length-$(n+1)$ polygons are globally joined. Therefore, a non-negligible, $n^{-o(1)}$, fraction of such polygons are charming at index~$k$.

By varying $k$ over order $n$ values,
we find that the snake method hypothesis~(\ref{e.afewcs}) is verified under the standing assumption, because it must be the case that, with probability at least $n^{-o(1)}$, the $\psap_{n+1}$-sampled polygon has at least $n^{1 - o(1)}$ indices lengths at which it is charming. Naturally, Theorem~\ref{t.snakesecondstep} then applies to show that the closing probability has a rapid decay, in contradiction to the standing assumption. 
  
In this way, the standing assumption is seen to be false, and the aim of showing that the closing probability decays as quickly as $n^{-1/2 - \eps}$ is achieved subsequentially (for some positive $\eps$).
The argument demonstrates both the beauty and horror of the author's experience of self-avoiding walk: it is rather delicate and intricate, and it achieves a specific inference, regarding closing probability, by examining in detail interesting geometric questions about walks and polygons,  concerning say the proportion of large polygons that are globally joined; however, it leaves these questions unsettled.

\section{Proving Theorem~\ref{t.thetexist}(2)}\label{s.theoremsecond}

The paper has seven further sections, and these are devoted to proving Theorem~\ref{t.thetexist}. We will gather the two set of preliminaries we have assembled and prove the theorem using the snake method implemented via polygon joining.

To summarise the structure of the rest of the paper, we mention first that we will begin with Theorem~\ref{t.thetexist}(2) because the proof of the first part of the theorem shares the framework of the second but is slightly more technical. 
We first reduce (in the rest of this section) Theorem~\ref{t.thetexist}(2)
to Proposition~\ref{p.thetexist}; in Section~\ref{s.threesteps}, we state the conclusions of the three main steps that lead to the proof of the proposition; and in the three ensuing sections, we prove the results associated to the respective steps and thereby complete the proof of Theorem~\ref{t.thetexist}(2). Section~\ref{s.prepare} presents a counterpart 
to Proposition~\ref{p.thetexist} for Theorem~\ref{t.thetexist}(1) and collects some useful statements for the proof of this part of the main theorem, which proof is presented with the same three step structure in the final Section~\ref{s.final}.

Focussing now on the proof of Theorem~\ref{t.thetexist}(2), note that it is enough to suppose the existence of  $\theta$ and $\xi$
and to argue that $\thet + \xi \geq 5/3$, since the latter assertion of  Theorem~\ref{t.thetexist}(2) was justified as it was stated.  In order to prove this inequality by finding a contradiction,  we suppose that $\thet + \xi < 5/3$.
The sought contradiction is exhibited in the next result.


\begin{proposition}\label{p.thetexist}
Let $d = 2$. Assume that the limits~$\thet : = \lim_{n \in 2\N} \thet_n$ and $\xi =  \lim_{n \in \N} \xi_n$ exist. Suppose further that $\thet + \xi < 5/3$.  Then, for some $c > 1$ and $\delta >0$, the set of $n \in 2\N + 1$ for which 
$$
\psaw_n \big( \Ga \closes \big) \leq c^{-n^{\delta}}
$$
intersects the shifted dyadic scale $\big[ 2^i - 1,2^{i+1} -1 \big]$ for all but finitely many $i \in \N$.
\end{proposition}
\noindent{\bf Proof of Theorem~\ref{t.thetexist}(2).}
The non-negative limits $\thet$ and $\xi$ are hypothesised to exist. If their sum is finite, then the formula  (\ref{e.closepc.formula}) exhibits a polynomial decay. 
When the sum is less than $5/3$, this 
  contradicts Proposition~\ref{p.thetexist}. \qed

It remains of course to derive the proposition. 
Assume throughout the derivation that the proposition's hypotheses are satisfied. Define an exponent $\chi$ so that  $\theta + \xi = 3/2 + \chi$. 
Since $\xi \geq 0$ and $\theta \geq 3/2$ (classically and by Theorem~\ref{t.polydev}), $\chi$ is non-negative.
Note that $\psaw_n \big( \Ga \closes \big)$ equals $n^{-1/2 - \chi + o(1)}$, and also that  we are supposing  that $\chi < 1/6$. 

The proof of Proposition~\ref{p.thetexist} is an application of the snake method.
We now set the snake method exponent parameters for this application, taking 
\begin{itemize}
\item the snake length $\beta$ equal to one; 
\item
the inverse charm $\alpha$ equal to $1/2 + 2\chi + 4\eps$;
\item and the deficit $\esm$ equal to $\chi + 8 \eps$;
\end{itemize}
 where henceforth in proving Proposition~\ref{p.thetexist}, $\eps > 0$ denotes a given but arbitrarily small quantity. The quantity $\delta = \beta - \esm - \alpha = 1/2 - 3\chi - 12\eps$ must be positive if the snake method is to work, and thus we choose $\eps < (1/6 - \chi)/4$.

\section{Three steps to the proof of Proposition~\ref{p.thetexist}}\label{s.threesteps}

The argument leading to Proposition~\ref{p.thetexist} has three principal steps. We explain these now in outline, stating the conclusion of each one. The proofs follow, in the three sections that follow this one.

\subsection{Step one: finding many first parts with high closing probability}

The outcome of the first step may be expressed as follows, using the high conditional closing probability set notation~$\hPhi$ introduced in Subsection~\ref{s.possparts}.
\begin{proposition}\label{p.notquite}
For $i \in \N$ sufficiently high, there exists $m' \in 2\N \cap \big[ 2^{i+3}, 2^{i+4} \big]$ such that, writing $\lset$ for the set of values $k \in \N$, $2^i \leq k \leq 2^i + 2^{i-2}$, that satisfy
$$
\psap_{m'+k} \Big( \Ga_{[0,k]} \in \hPhi_{k,k+m' - 1}^\alpha \Big)
 \geq 1 - 2^{-i\chi} \, ,
$$
we have that $\vert \lset \vert \geq 2^{-8}m'$.
\end{proposition}

This first step is a soft, Fubini, argument: such index pairs $(k,m'+k)$ are characteristic given our assumption on closing probability decay, and the proposition will follow from Lemma~\ref{l.closecard}.

\subsection{Step two: individual snake terms are often charming}
To apply the snake method, we must verify its fundamental hypothesis~(\ref{e.afewcs}) that charming snakes are not rare under the uniform polygon law. Proposition~\ref{p.notquite}
is not adequate for verifying this hypothesis, because the polygon law index $m' + k$ is a variable that changes with $k$.
We want a version of the proposition where this index is fixed. 
The next assertion, which is the conclusion of our second step, is such a result. Note from the differing forms of the right-hand sides in the two results 
that a behaviour determined to be highly typical in Proposition~\ref{p.notquite} has a counterpart in Proposition~\ref{p.almost} which is merely shown to be not unusual.
\begin{proposition}\label{p.almost}
For each $i \in \N$ sufficiently high, there exist $m \in 2\N \cap \big[ 2^{i+4}, 2^{i+5} \big]$ and $m' \in \big[ 2^{i+3},2^{i+4} \big]$ such that, writing $K$ for the set of values $k \in \N$, $1 \leq k \leq m-m'$, that satisfy
$$
\psap_m \Big( \Ga_{[0,k]} \in \hPhi_{k,k+m' - 1}^\alpha \Big)
 \geq   m^{-\esm + \eps} \, ,
$$
we have that $\vert K \vert \geq 2^{-9} m$.
\end{proposition}

\begin{figure}
  \begin{center}
    \includegraphics[width=1.0\textwidth]{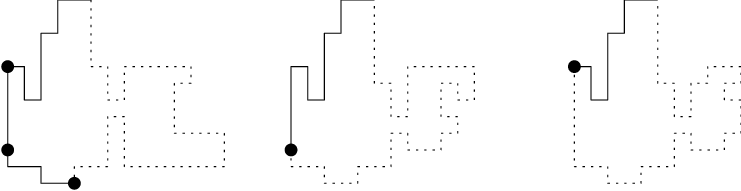}
  \end{center}
  \caption{The snake method via polygon joining. Black dots mark vertices of  polygons $\ga$: fixing a typical $\tilden \in 2\N$, assign a black dot at each vertex $\ga_j$, $\tilden/4 \leq j \leq 3\tilden/4$, of a polygon $\ga$ of length of order~$\tilden$, that has the property that $\saw^0_{\tilden+j} \big( \Ga \closes \, \big\vert \, \vert \Ga^1 \vert = j , \Ga^1 = \ga_{[0,j]} \big) \geq (\tilden)^{-1/2 - \chi - \eps}$.
Our assumptions readily imply that most polygons $\ga$ with length $n'$ of order say~$2\tilden$ will have a black dot at $\ga_{n' - \tilden}$. What about at other locations $\ga_j$ where $j$ is also of order $\tilden$? A black dot is likely to appear at $\Ga_j$ when $\Ga$ is $\psap_{j+\tilden}$-distributed. If we can make a {\em comparison of measure}, showing that the laws $\psap_{j+\tilden}$ and $\psap_{{n'}}$ are to some degree similar, then the black dot will be known to also appear at~$\ga_j$ in a typical sample of $\psap_{{n'}}$. Applying this for many such~$j$, black dots will appear along the course of a typical length~$n'$ polygon. These black dots index charming snake terms and permit the use of the snake method.
Comparison of measure  will be undertaken in Proposition~\ref{p.similarity} by the use of polygon joining. Lemma~\ref{l.polyepsregthet} shows that a non-negligible proportion of polygons are regulation global join polygons; thus, such polygons themselves typically have black dots~$\tilden$ steps from their end.
The law of the left polygon under the uniform law on regulation polygons of a given length is largely unchanged as that total length is varied on a given dyadic scale,
because the length discrepancy can be absorbed by altering the length of the right polygon. Indeed, in the above sketches, we
deplete the length of the right polygon in a possible extension of the depicted first part, as the length of this first part shortens; in this way, we show that first parts typically arising in regulation polygons at one length are also characteristic in such polygons with lengths on the same scale.}\label{f.snakemethodpolygon}
\end{figure}

The second step, leading to the proof of this proposition, is quite subtle. Expressed in its simplest terms, the idea of the proof is  that the subscript index in $\psap$ can be changed from the $k$-dependent $m' + k$ in Proposition~\ref{p.notquite}  to the $k$-independent~$m$ in the new result by using the polygon joining technique.
This index change involves proving a {\em similarity of measure} between polygon laws with distinct indices on a given dyadic scale. The set of regulation global join polygons introduced in Section~\ref{s.pjprep} provides a convenient reservoir of macroscopically joined polygons whose left polygons are shared between polygon laws of differing index, so that this set provides a means by which this similarity may be proved. 
Figure~\ref{f.snakemethodpolygon} offers some further outline of how we will exploit the ample supply of regulation polygons to prove similarity of measure.

\subsection{Step three: applying the snake method Theorem~\ref{t.snakesecondstep}}
In this step, it is our goal to use Proposition~\ref{p.almost}
in order to verify the snake method hypothesis~(\ref{e.afewcs}). 
In the step, we will  fix the method's two index parameters: for $i \in \N$ given, we will use Proposition~\ref{p.almost} to take $n$ equal to~$m-1$, and $\ell$ equal to $m - m'$. 

We may now emphasise why the proposition is useful for this verification goal.
Note that $\max K \leq \ell$.
Since $n - \ell = m' - 1$, we have that, for $k \in [1,\ell]$, the walk $\gamma \in \fpart_{\ell,n}$ is $(\alpha,n,\ell)$-charming at index $k$ if and only if 
$$
 \psaw^0_{k+m'-1} \Big( \Ga \closes \, \Big\vert \, \big\vert \Gamma^1 \big\vert = k \, , \,  \Ga^1_{[0,k]} = \ga_{[0,k]} \Big) > n^{-\alpha} \, ;
$$ 
note the displayed condition is almost the same as $\ga_{[0,k]} \in \hPhi_{k,m' +k-1}^{\alpha}$, the latter occurring when $n^{-\alpha}$ is replaced by~$(m'+k-1)^{-\alpha}$, a change which involves only a bounded factor. 

Thus, Proposition~\ref{p.almost} shows that individual snake terms are charming with probability at least~$m^{-\chi - o(1)}$.
It will be a simple Fubini argument that will allow us to verify~(\ref{e.afewcs}) by confirming that it is not rare that such terms gather together to form charming snakes. Theorem~\ref{t.snakesecondstep} may be invoked to prove Proposition~\ref{p.thetexist} immediately. The third step will provide this Fubini argument.

The next three sections are devoted to implementing these three steps. In referring to the dyadic scale index $i \in \N$ in results stated in these sections, we will sometimes neglect to record that this index is assumed to be sufficiently high.

\section{Step one: the proof of Proposition~\ref{p.notquite}}\label{s.notquite}

Before proving this result, we review its meaning in light of  Figure~\ref{f.snakemethodpolygon}. 
Proposition~\ref{p.notquite} is a precisely stated counterpart to the assertion in the figure's caption that a black dot typically appears $\tilden$ steps from the end of a polygon drawn from $\psap_{n'}$ where $n'$ has order $2\tilden$ (or say $2^{i+3}$). Note the discrepancy that we have chosen $\alpha$ to be marginally above $1/2 + 2\chi$ while in the caption a choice just above $1/2 + \chi$ was made. The extra margin permits $1 - 2^{-i\chi}$ rather than $1 - o(1)$ in the conclusion of the proposition. We certainly need the extra margin; we will see how  at the end of step two.

Consider two parameters $\eclo > \eroom > 0$. (During the proof of Proposition~\ref{p.almost}, we will ultimately set  $\eclo = 4\eroom$ and $\eroom = \eps$, where $\eps > 0$ is the parameter that has been used to specify the three snake method exponents.)

\begin{definition}\label{d.cpt}
Let $i \in \N$. An index pair $(k,j) \in \N \times 2\N$ will be called {\em closing probability typical} on dyadic scale $i$ if 
\begin{itemize}
\item $k \in \big[ 2^i , 2^i + 2^{i-2} \big]$ and $j \in \big[ 2^{i+4},  2^{i+5} \big]$;
\item and, for any $a \geq 0$,
\begin{equation}\label{e.acondition}
 \psap_{j} \Big( \Gamma_{[0,k]} \notin \phicl{k}{j-1}{1/2 + (a+1)\chi + \eclo}  \Big) \leq  (j-1)^{-(a\chi + \eclo - \eroom)} \, .
\end{equation}
\end{itemize}
\end{definition}

\begin{proposition}\label{p.manyethet}
For all $i \in \N$ high enough, there exists $k \in 2\N \cap \big[ 2^{i+4}, 2^{i+4} + 2^i \big]$ for which there are at least $2^{i-4}$ values of $j \in 2\N \cap [0,2^{i-2}]$ with the index pair
$\big( 2^i + j, k + j \big)$ closing probability typical on dyadic scale $i$.
\end{proposition}

When we prove Proposition~\ref{p.almost} in step two, we will in fact do so using a direct consequence of Proposition~\ref{p.manyethet}, rather than Proposition~\ref{p.notquite}. 
The latter result, a byproduct of the proof of  Proposition~\ref{p.manyethet}, has been stated merely because it permitted us to make a direct comparison between steps one and two in the preceding outline.
As such, our goal in step one is now to prove Proposition~\ref{p.manyethet}.
\begin{definition}\label{d.e}
Define $E$ to be the set of pairs $(k,j) \in \N \times 2\N$,  
\begin{itemize}
\item where $k \in \big[ 2^i , 2^i + 2^{i-2} \big]$ and $j \in \big[ 2^{i+4},  2^{i+5} \big]$;
\item and the pair $(k,j)$ is such that
\begin{eqnarray}
 & & \# \Big\{ \gamma \in \saw^0_{j-1}: \vert \gamma^1 \vert = k \Big\} \label{e.rchithet} \\
 & < & (j-1)^{1/2 + \chi + \eroom} \cdot  \# \Big\{ \gamma \in \saw^0_{j-1}: \vert \gamma^1 \vert = k , \gamma \closes \Big\} \, . \nonumber
\end{eqnarray}
\end{itemize}
\end{definition}

We will now apply Lemma~\ref{l.closecard} to show that membership of $E$ is typical. The application is possible because by hypothesis there is a lower bound on the  decay of the closing probability. Indeed, we know that
$\psaw_n \big( \Gamma \closes \big) \geq n^{-1/2 - \chi - o(1)}$. Thus, the lemma with $\alphamac = 1/2 + \chi + \eroom/2$ and $\deltamac = \eroom/2$ proves the following.
\begin{lemma}\label{l.etypical}
Provided the index $i \in \N$ is sufficiently high,
for each $j \in  2 \N \cap \big[ 2^{i+4},  2^{i+5} \big]$, the set of $k \in \big[ 2^i , 2^i + 2^{i-2} \big]$ such that $(k,j) \not\in E$ has cardinality at most $2 j^{1 - \eroom/2} \leq 2 \big(  2^{i+5} \big)^{1-\eroom/2}$.
\end{lemma}

\begin{lemma}\label{l.rjsjthet}
Any pair $(k,j) \in E$
is closing probability typical on dyadic scale $i$.
\end{lemma}
\noindent{\bf Proof.} We must verify that $(k,j)$
verifies~(\ref{e.acondition}) for any $a \geq 0$. Recalling the notation $\psaw_{j-1}^0$ from Subsection~\ref{s.sawfree}, note that the assertion
\begin{equation}\label{e.lststhet}
 \psaw^0_{j-1} \Big( \Gamma^1 \notin \phicl{k}{j-1}{1/2 + (a+1)\chi + \eclo} \, \, \Big\vert \, \, \Gamma \closes \, ,  \vert \Gamma^1 \vert  = k \Big) \leq (j - 1)^{-(a\chi + \eclo - \eroom)} 
\end{equation}
is a reformulation of this condition.
Note that
\begin{equation}\label{e.cardgamma}
  \# \Big\{ \gamma \in \saw^0_{j-1}: \vert \gamma^1 \vert =  k \, , \, \gamma^1 \notin \phicl{k}{j-1}{1/2 + (a+1)\chi + \eclo} \Big\} 
\end{equation}
is at least the product of $(j-1)^{1/2 + (a+1)\chi + \eclo}$
and
$$
 \# \Big\{ \gamma \in \saw^0_{j-1}: \vert \gamma^1 \vert = k \, , \, \gamma^1 \notin \phicl{k}{j-1}{1/2 + (a+1)\chi + \eclo} \, , \, \gamma \closes \Big\}  \, .
$$
The quantity~(\ref{e.cardgamma}) is bounded above by~(\ref{e.rchithet}); thus is~(\ref{e.lststhet}) obtained. \qed

\medskip

\noindent{\bf Proof of Proposition~\ref{p.manyethet}.}
By Lemma~\ref{l.rjsjthet}, it is enough to argue that there exists $k \in 2\N \cap \big[ 2^{i+4}, 2^{i+4} + 2^i \big]$ for which there are at least $2^{i-4}$ values of $j \in 2\N \cap [0,2^{i-2}]$ with 
$\big( 2^i + j, k + j \big) \in E$.
 
By Lemma~\ref{l.etypical},
\begin{eqnarray}
 & & \sum_{j \in 2\N \cap [ 2^{i+4} + 2^{i-2} , 2^{i+4} + 2^{i} ] } \sum_{\ell = 2^i}^{2^i + 2^{i-2}} 
 1\!\!1_{(\ell,j) \in E} \label{e.ineqe} \\
 & \geq & \Big( 2^{i-2} + 1 - 2 \big( 2^{i+5} \big)^{1 - \eroom/2} \Big) \cdot \tfrac{1}{2} \big( 2^{i} - 2^{i-2} \big) \, . \nonumber
\end{eqnarray}
    
\begin{figure}
    \begin{center}
      \includegraphics[width=0.5\textwidth]{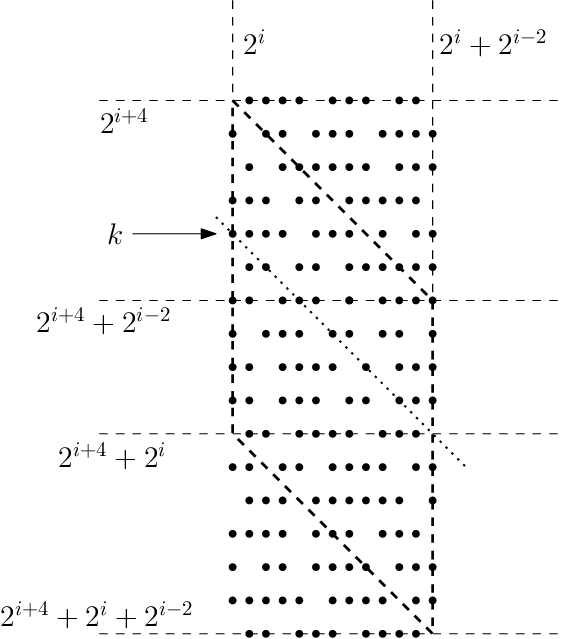}
    \end{center}
    \caption{Illustrating the proof of Proposition~\ref{p.manyethet}.    The dots depict elements of $E$. 
 The quantity~(\ref{e.ineqe})  is the number of dots in the middle rectangle; (in these terms, $\ell$ is the $x$-coordinate in the sketch and $j$ is the negative $y$-coordinate). Similarly,     
    the quantity~(\ref{e.doublesum}) is the number of dots in the parallelogram bounded by dashed lines.
    In~(\ref{e.kterm}), we find a diagonal, such as the one depicted with a pointed line, on which many dots lie. }\label{f.sumexplain}
  \end{figure}

As Figure~\ref{f.sumexplain} illustrates, the left-hand side here is bounded above by 
\begin{equation}\label{e.doublesum}
 \sum_{k \in 2\N \cap [ 2^{i+4},2^{i+4} + 2^{i} ]} \sum_{\ell = 2^i}^{2^i + 2^{i-2}}
 1\!\!1_{(\ell,k + \ell - 2^i) \in E} \, .
\end{equation}
There being $2^{i - 1} + 1 \leq 2^{i}$ indices $k \in 2\N \cap \big[ 2^{i+4},2^{i+4} + 2^{i} \big]$, one such~$k$ satisfies
\begin{equation}\label{e.kterm}
\sum_{\ell = 2^i}^{2^i + 2^{i-2}} 1\!\!1_{(\ell,k + \ell - 2^i) \in E} \geq \tfrac{3}{8} \, \Big( 2^{i-2} + 1 - 2 \big( 2^{i+5} \big)^{1 - \eroom/2} \Big) \, .
\end{equation}
Noting that $\tfrac{3}{4} \big( 2^{i+5} \big)^{1 - \eroom/2} \leq \tfrac{1}{8} 2^{i-2}$, we obtain the sought statement and so conclude the proof. \qed

\medskip

\noindent{\bf Proof of Proposition~\ref{p.notquite}.} 
Proposition~\ref{p.manyethet} furnishes the existence of a minimal $Q \in 2\N \cap \big[ 2^{i+4}, 2^{i+4} + 2^{i} \big]$ 
such that there are at least $2^{i-4}$ values of $k \in 2\N \cap [0,2^{i-2}]$ with  $\big( 2^i + k, Q + k \big)$ closing probability typical on dyadic scale~$i$. 

Let $\big( k_1,\ldots,k_{2^{i-4}} \big)$ be an increasing sequence of such values of $k$ associated to the value~$Q$. 
Set $m' = Q - 2^i$ and note that $m' \in \big[2^{i+3},2^{i+4} \big]$. 
Take $a=1$ in~(\ref{e.acondition}) and specify that $\eclo  = 4\eps$ alongside $\eclo > \eroom$  to find that, for all $j \in [1,2^{i-4}]$,  $2^i + k_j \in \lset$. Thus, $\vert \lset \vert \geq 2^{i-4} \geq 2^{-8} m'$.

This completes the proof of the proposition. \qed

\medskip

We end this section by recording the explicit artefact that will be used in step two.
\begin{definition}\label{d.rjsj}
Set $\indexc = 2^{i-4}$. We specify a {\em constant difference} sequence of  index pairs $(r_j,s_j)$, $1 \leq j \leq \indexc$. Its elements are closing probability typical for dyadic scale $i$, the difference $s_j - r_j$ is independent of $j$ and takes a value in $[2^{i+4}- 2^i,2^{i+4}]$, and $\big\{ r_j: 1 \leq j \leq \indexc \big\}$ is increasing.
Moreover, $r_j \in 2\N \cap [2^i,2^i + 2^{i-2}]$ and  $s_j \in 2\N \cap [ 2^{i+4} , 2^{i+4} + 2^i + 2^{i-2}  ]$ for each $j \in [1,\indexc]$.
\end{definition}
Indeed, we may work with the sequence constructed in the preceding proof, 
setting $r_j = 2^i + k_j$ and $s_j = Q + k_j$  for $j \in \big[ 1, \indexc  \big]$. 

We are {\em en route} to Proposition~\ref{p.almost},
and take the opportunity to mention how its two parameters $m$ and $m'$ will be specified. The latter will be chosen, as it was in the proof of Proposition~\ref{p.notquite}, to be the value of the constant difference in the sequence in Definition~\ref{d.rjsj}, while the former will be taken equal to $s_{\indexc}$ in this definition.

\section{Step two: deriving Proposition~\ref{p.almost}}\label{s.almost}

The key tool for this proof is similarity of measure
between polygon laws with distinct indices, with the reservoir of regulation global join polygons providing the mechanism for establishing this similarity. Recall that these polygons were introduced and discussed in Section~\ref{s.pjprep}; the reader may wish to review this material~now. 

This step two section has three subsections. In the first, we provide a piece of apparatus, which augments the regulation polygon properties discussed in Section~\ref{s.pjprep}, that we will need in order to state our similarity of measure result. This result, Proposition~\ref{p.similarity}, is stated and proved in the next subsection. The final subsection gives the proof of Proposition~\ref{p.almost}.

\subsection{Relating polygon laws $\psap_n$ 
for distinct $n$ via global join polygons}\label{s.regnotrare}

 The next lemma shows that there is an ample supply of regulation polygons in the sense that they are not rare among all polygons.
In contrast to the regulation polygon tools in Section~\ref{s.pjprep}, Lemma~\ref{l.polyepsregthet}
is contingent on our standing assumption that the hypotheses of Proposition~\ref{p.thetexist} are satisfied.
(In fact, it is the existence and finiteness of $\theta$ and $\xi$ that are needed, rather than the condition $\theta + \xi < 5/3$.)

Define $\reg_{i} : = 2\N \cap \big[ 2^{i}, 2^{i+1} \big]$. 

\begin{lemma}\label{l.polyepsregthet}
Let $\varphi > 0$. 
For any $i \in \N$ sufficiently high, and $n \in 2\N \cap \big[ 2^{i+4},2^{i+5} \big]$,
$$
\psap_n \bigg( \Gamma \in \bigcup_{k \in \reg_{i+2}} \dubbub_{k,n-16-k}  \bigg) \geq n^{-\chi - \varphi} \, .
$$
\end{lemma}

 The set  $\reg_{i}$ may be thought as a set of {\em regular} indices $j \in 2\N \cap \big[ 2^{i}, 2^{i+1} \big]$, with~$\thet_j$ for example being neither atypically high nor low.  
As we prove Theorem~\ref{t.thetexist}(2), all indices are regarded as regular, so that we make the trivial specification  $\reg_{i}  = 2\N \cap \big[ 2^{i}, 2^{i+1} \big]$; the set will be respecified non-trivially in the proof of Theorem~\ref{t.thetexist}(1).


(As we explained in the three-step outline, the snake method parameters, including $n$, will be set for the proof of Proposition~\ref{p.thetexist} at the start of step three. Until then, we are treating $n$ as a free variable.)

\medskip

\noindent{\bf Proof of Lemma~\ref{l.polyepsregthet}.} Note that the probability in question is equal to the ratio of the cardinality of the union $\bigcup_{j \in \reg_{i+2}} \dubbub_{j,n-16-j}$ and the polygon number $p_n$.

In order to bound below the numerator in this ratio, we aim to apply Proposition~\ref{p.polyjoin.aboveonehalf} with the role of $i \in \N$ played by $i+2$ and with $\mathsf{R} = \reg_{i+2}$.
In order to do so, set $\overline{\theta} = \sup_{n \in 2\N} \theta_n$, a quantity that is finite by hypothesis. 
Set the proposition's $\Theta$ equal to  $\overline{\theta}$. 
The index $\kmac$ in the proposition may then be any element of $\reg_{i+2}$.
Applying the proposition, we find that there exists $c = c(\Theta) > 0$ such that
$$
 \Big\vert \, \bigcup_{j \in \reg_{i+2}} \dubbub_{j,n-16-j} \, \Big\vert \, \geq \, c \,
      \frac{n^{1/2}}{\log n} \,
\sum_{j \in \reg_{i+2}} p_j p_{n-16-j} \, .
$$

The existence of~$\thet$
implies that $p_m \mu^{-m} \in \big[ m^{-\thet - \varphi} , m^{-\thet + \varphi} \big]$ whenever $m \geq m_0(\varphi)$.
Since $i \in \N$ is supposed sufficiently high, we thus find that
\begin{eqnarray*}
& & \psap_n \Big(  \Gamma \in \bigcup_{j \in \reg_{i+2}} \dubbub_{j,n-16-j}   \Big) \\
 & \geq & \mu^{-16} c \, \frac{n^{1/2}}{\log n}  \sum_{j \in \reg_{i+2}} j^{-(\thet + \varphi)} (n-16-j)^{-(\thet + \varphi)} \cdot n^{\theta - \varphi} \geq 2^{-4} \mu^{-16}  c \, \frac{n^{3/2 - \thet - 3\varphi}}{\log n} \, ,
\end{eqnarray*}
where we used $\# \reg_{i+2} \geq 2^{i+1} \geq 2^{-4}n$. Noting that $\thet \leq 3/2 + \chi$ and relabelling~$\varphi > 0$ completes the proof. \qed

\subsection{Similarity of measure}


We may now state our similarity of measure Proposition~\ref{p.similarity}. Recall from Definition~\ref{d.rjsj} the constant difference sequence $(r_j,s_j)$, $1 \leq j \leq \indexc$, and  also the regular index set $\reg_{i+2}$ in the preceding subsection.
\begin{proposition}\label{p.similarity}
For each $j \in [1,L]$, there is a subset $D_j \subseteq \fpart_{r_j}$ 
satisfying
\begin{equation}\label{e.similarity}
 \psap_{s_j} \Big( \Gamma_{[0,r_j]} \not\in D_j \, , \, \Gamma \in \dubbub_{k,s_j - 16 - k} \, \, \, \textrm{for some $k \in \reg_{i+2}$} \Big) \leq 4 \mu^{16}  s_j^{-10} 
\end{equation}
such that, for each $\phi \in D_j$,
\begin{eqnarray*}
 & & \psap_{s_{\indexc}} \Big( \Gamma_{[0,r_j]} = \phi  \Big) \\
 & \geq & \tfrac{1}{40} \,  \conindex^{-1} \big( \log s_\indexc \big)^{-1} \big( \tfrac{3}{16} \big)^{\thet} \, s_{\indexc}^{-4\eps}  \\
 & & \qquad \qquad \qquad \times \, \, \, 
  \psap_{s_j} \Big(  \Gamma_{[0,r_j]} = \phi \, , \,  \Gamma \in \dubbub_{k,s_j - 16 - k} \, \, \textrm{for some $k \in \reg_{i+2}$} \Big) \, .
\end{eqnarray*}
\end{proposition}

Preparing for the proof, we define, for each $\gamma \in \sap_{s_j}$, $1 \leq j \leq \indexc$, the set of $\ga$'s regulation global join indices 
$$
\tau_\gamma = \Big\{ k \in \reg_{i+2} : \gamma \in \dubbub_{k,s_{\indexcj} - 16 - k} \Big\} \, .
$$

 For $\gamma \in \sap_{s_{\indexcj}}$, the map sending $k \in \tau_\gamma$ to the junction plaquette associated to the join of elements of $\sapell_k$ and $\sapr_{s_\indexcj - 16 - k}$
is an injective map into $\globaljoin_\ga$: see the first assertion of Proposition~\ref{p.rgjpnumber}. Thus,  
$\vert \tau_\gamma \vert \leq \big\vert \globaljoin_\gamma \big\vert$.  Corollary~\ref{c.globaljoin} implies that   there exists a constant $\conindex > 0$ such that
\begin{equation}\label{e.constl}
 \psap_{s_{\indexc}} \Big( \big\vert \tau_\Gamma \big\vert \geq \conindex \log s_{\indexc} \Big) \leq \conindex \, s_{\indexc}^{-4\overline\theta - 22} \, ,
\end{equation}
where here we have reintroduced $\overline\theta = \sup_{n \in 2\N} \theta_n$, a quantity which we recall is finite by hypothesis; (this finiteness permits the application of the corollary).

We make use of this constant to specify for each $j \in [1,\indexc]$ a
set $D_j \subseteq \fpart_{r_j}$ of length~$r_j$ first parts that, when extended to form a length~$s_\indexc$ polygon, do not typically produce far-above-average numbers of regulation joinings. We define
$$
D_j = \Big\{ \phi \in \fpart_{r_j}: \psap_{s_{\indexc}} \Big( \big\vert \tau_\Gamma \big\vert \leq \conindex \log s_{\indexc} \, \Big\vert \, \Gamma_{[0,r_j]} = \phi    \Big) \geq 1 - \conindex \, s_{\indexc}^{-2 \overline\theta - 11} \Big\} \, ,
$$

The use of the next lemma lies in its final part, which establishes the $\psap_{s_j}$-rarity
assertion~(\ref{e.similarity}).
The first two parts of the lemma are steps towards the third.
The rarity under  $\psap_{s_\indexc}$ of the length $r_j$ initial path failing to be in $D_j$  is a direct consequence of definitions, as the first part of the lemma shows, but a little work is needed to make the change of index $\indexc \to j$.

\begin{lemma}\label{l.dcub}
Let $j \in [1,\indexc]$. 
\begin{enumerate}
\item We have that 
$$
\psap_{s_{\indexc}} \Big( \Gamma_{[0,r_j]} \in D_j^c \Big) \leq  s_{\indexc}^{-2 \overline\theta - 11} \, .
$$
\item
Let $k \in \reg_{i+2}$. Then
$$
\psapleft k \Big( \Gamma_{[0,r_j]} \not\in D_j \Big)
\leq   16 \mu^{16}  s_\indexc^{-11} \, .
$$
\item We have that
$$
 \psap_{s_j} \Big( \Gamma_{[0,r_j]} \not\in D_j \, , \, \Gamma \in \dubbub_{k,s_j - 16 - k} \, \, \, \textrm{for some $k \in \reg_{i+2}$} \Big) \leq 4 \mu^{16}  s_j^{- 10} \, .
$$
\end{enumerate}
\end{lemma}
\noindent{\bf Proof (1).}
Note that
\begin{eqnarray*}
 & & \psap_{s_{\indexc}} \Big( \big\vert \tau_\Gamma \big\vert \geq \conindex \log s_{\indexc} \Big) \geq \psap_{s_{\indexc}} \Big( \big\vert \tau_\Gamma \big\vert \geq \conindex \log s_{\indexc} \, , \, \Gamma_{[0,r_j]} \in D_j^c \Big) \\
  & =  &   \psap_{s_{\indexc}} \Big(\Gamma_{[0,r_j]} \in D_j^c \Big)  \psap_{s_{\indexc}} \Big( \big\vert \tau_\Gamma \big\vert \geq \conindex \log s_{\indexc} \, \Big\vert \, \Gamma_{[0,r_j]} \in D_j^c \Big) \\
 &  \geq & \conindex s_{\indexc}^{-2 \overline\theta - 11} \,  \psap_{s_{\indexc}} \Big(\Gamma_{[0,r_j]} \in D_j^c \Big)  \, . 
\end{eqnarray*}
Lemma~\ref{l.dcub}(1) thus follows from (\ref{e.constl}).  

\medskip

\noindent{\bf (2).}
In light of the preceding part, it suffices to prove 
\begin{equation}\label{e.sjsl}
\psapleft{k} \Big( \Gamma_{[0,r_j]} \not\in D_j \Big)
\leq  16 \mu^{16} s_\indexc^{2\overline{\theta}} \,
\psap_{s_\indexc} \Big( \Gamma_{[0,r_j]} \not\in D_j \Big) \, .
\end{equation}
We will derive this bound with a simple joining argument. 
Note that $s_\indexc \geq k + 16$ since $s_\indexc \geq 2^{i+4}$, $k \leq 2^{i+3}$ and $i \geq 1$.
Consider the map 
$$
\Psi:  \big\{ \phi \in \sapell_k : \phi_{[0,r_j]} \in D_j^c \big\} \times  \sapr_{s_\indexc - k - 16} \to  \big\{ \phi \in \sap_{s_\indexc} : \phi_{[0,r_j]} \in D_j^c \big\}
$$
defined by setting $\Psi(\phi_1,\phi_2)$ equal to the join polygon $J(\phi_1,\phi_2 + \vec{u})$, where the second argument is chosen so that the pair is globally Madras joinable. (A choice of $\vec{u}$ ensuring that the pair is {\em globally} joinable may be made by placing  the second polygon strictly below the $x$-axis and then displacing it horizontally to the rightmost location at which it is Madras joinable to the first polygon.)   

We now wish to apply Lemma~\ref{l.leftlongbasic} 
to learn that 
   $\Psi ( \phi , \phi' )_{[0,r_j]} =  \phi_{[0,r_j]}$, so that $\Psi$'s image lies in the stated set. In order to do so, it is enough to verify that $r_j \leq k/2 - 1$, since $\phi \in \sap_k$ is left-long. This bound follows from $k \geq 2^{i+2}$, $r_j \leq 2^i + 2^{i-2}$
   (as noted in Definition~\ref{d.rjsj}), and $i \geq 1$.

   The decomposition uniqueness in Proposition~\ref{p.rgjpnumber} shows that $\Psi$ is injective.
   Thus, by considering $\Psi$, we find that
$$
\# \big\{ \phi \in \sap_{s_\indexc} : \phi_{[0,r_j]} \in D_j^c \big\} \geq 
\# \big\{ \phi \in \sapell_k : \phi_{[0,r_j]} \in D_j^c \big\}  \# \big\{ \phi \in \sapr_{s_\indexc - k - 16} \big\} \, .  
$$   
By Lemma~\ref{l.polysetbound},
$$
 \psap_{s_\indexc} \big( \Gamma_{[0,r_j]} \in D_j^c \big)
  \geq \tfrac{1}{16} \psapleft{k} \big( \Gamma_{[0,r_j]} \in D_j^c \big) p_{s_\indexc}^{-1} p_k p_{s_\indexc - k - 16} \, .
$$  
Note that $p_{s_\indexc}^{-1} p_k p_{s_\indexc - k - 16} \geq \mu^{-16} s_\indexc^{-2\overline{\thet}}$, where recall that $\overline{\thet} = \sup_{n \in 2\N} \thet_n$. In this way, we obtain~(\ref{e.sjsl}). 

\medskip

\noindent{\bf (3).} 
 We first make the key observation that,
for $k \in \reg_{i+2}$ and $j' \in [1,\indexc]$,
\begin{equation}\label{e.keyob}
 \psap_{s_{j'}} \Big( \Gamma_{[0,r_j]} = \phi \, \Big\vert \, \Gamma \in \dubbub_{k,s_{j'} - 16 - k} \Big)  =  \psapleft{k} \Big( \Gamma_{[0,r_j]} = \phi \Big) \, .
\end{equation}
Proposition~\ref{p.leftlong} implies this, provided that its hypotheses that $r_j \leq k/2 - 1$ and $k/2 \leq s_{j'} - 16 - k \leq 35 k$ are valid. 
That the first inequality holds has been noted in the proof of the lemma's second part. That the second holds follows from $2^{i+3} \geq k \geq 2^{i+2}$, $2^{i+5} \geq s_{j'} \geq 2^{i+4}$ and $i \geq 2$.

The quantity on the left-hand side of the inequality in Lemma~\ref{l.dcub}(3) is at most
\begin{eqnarray*}
 & & \sum_{k \in \reg_{i+2}} \psap_{s_j} \Big( \Gamma_{[0,r_j]} \not\in D_j \, , \, \Gamma \in \dubbub_{k,s_j - 16 - k}  \Big)
\end{eqnarray*}
and thus by~(\ref{e.keyob}) with $j' = j$, at most
$\vert \reg_{i+2} \vert \max_{k \in \reg_{i+2}} 
\psapleft{k} \big( \Gamma_{[0,r_j]} \not\in D_j \big)$. 
We may now apply the second part of the lemma as well as the bounds
$\vert \reg_{i+2} \vert \leq 4^{-1} s_\indexc$ and $s_\indexc \geq s_j$ to obtain Lemma~\ref{l.dcub}(3). 
(This upper bound on $\vert \reg_{i+2} \vert$ is due to $\vert \reg_{i+2} \vert = 2^{i+1} + 1 \leq 2^{i+2}$ and $s_\indexc \geq 2^{i+4}$.) 
\qed

\begin{lemma}\label{l.compare}
Recall that $\indexc = 2^{i-4}$. 
Let $j \in [1,\indexc]$. 
\begin{enumerate}
\item 
For $\phi \in D_j$,
$$   
   \psap_{s_{\indexc}} \Big( \Gamma_{[0,r_j]} = \phi  \Big) \\
  \geq 
      \tfrac{1}{400} \, \conindex^{-1} \big( \log s_\indexc \big)^{-1} s_{\indexc}^{-2\eps} \mu^{-16} \sum_{k \in \reg_{i+2}} \psapleft{k} \Big( \Gamma_{[0,r_j]} = \phi  \Big)
   k^{1/2} p_k \mu^{-k} \, .
   $$ 
\item For $\phi \in \fpart_{r_j}$,
\begin{eqnarray*}
   & & \psap_{s_j} \Big(  \Gamma_{[0,r_j]} = \phi \, , \,  \Gamma \in \dubbub_{k,s_j - 16 - k} \, \, \textrm{for some $k \in \reg_{i+2}$} \Big) \\
   & \leq &  \tfrac{1}{10} \big( \tfrac{16}{3} \big)^{\thet} \mu^{-16} \cdot s_j^{2\eps}  \sum_{k \in \reg_{i+2}} 
   \psapleft{k}  \Big(  \Gamma_{[0,r_j]} = \phi \Big)  k^{1/2} p_k \mu^{-k} \, .  
\end{eqnarray*}
\end{enumerate}
\end{lemma}
\noindent{\bf Proof:} {\em (1).}
Note that $\psap_{s_{\indexc}} \big( \Gamma_{[0,r_j]} = \phi  \big)$ is at least
\begin{eqnarray*}
 & &  \psap_{s_{\indexc}} \Big( \Gamma_{[0,r_j]} = \phi \, , \, \Gamma \in \bigcup_{k \in  \reg_{i+2}} \dubbub_{k,s_{\indexc} - 16 - k}  \, , \, \big\vert \tau_{\Gamma} \big\vert \leq \conindex \log s_\indexc \Big) \\
 & \geq & \big( \conindex \log s_\indexc \big)^{-1}
 \sum_{k \in  \reg_{i+2}}
 \psap_{s_{\indexc}} \Big( \Gamma_{[0,r_j]} = \phi \, , \, \Gamma \in \dubbub_{k,s_{\indexc} - 16 - k}  \, , \, \big\vert \tau_{\Gamma} \big\vert \leq \conindex \log s_\indexc \Big)  
  \, . 
\end{eqnarray*}
The summand in the last line may be written
\begin{eqnarray*}
 & & \psap_{s_{\indexc}} \Big( \Gamma_{[0,r_j]} = \phi \, , \, \Gamma \in \dubbub_{k,s_{\indexc} - 16 - k}  \Big)  \\
  & - &
 \psap_{s_{\indexc}} \Big( \Gamma_{[0,r_j]} = \phi \, , \, \Gamma \in \dubbub_{k,s_{\indexc} - 16 - k}  \, , \, \big\vert \tau_{\Gamma} \big\vert > \conindex \log s_\indexc \Big)  \, ,
\end{eqnarray*}
with the subtracted term here being at most
$$
  \psap_{s_{\indexc}} \Big( \Gamma_{[0,r_j]} = \phi \, , \, \big\vert \tau_{\Gamma} \big\vert > \conindex \log s_\indexc \Big) \leq   \conindex s_{\indexc}^{-2 \overline\theta - 11}  \psap_{s_{\indexc}} \Big( \Gamma_{[0,r_j]} = \phi \Big) 
$$
since $\phi \in D_j$.   
We noted at the end of the preceding proof that $\vert \reg_{i+2} \vert \leq 4^{-1} s_\indexc$.
Thus,  $\psap_{s_{\indexc}} \big( \Gamma_{[0,r_j]} = \phi  \big)$ is at least
\begin{eqnarray}
 & & \big( \conindex \log s_\indexc \big)^{-1}
 \sum_{k \in  \reg_{i+2}} \psap_{s_{\indexc}} \Big( \Gamma_{[0,r_j]} = \phi \, , \, \Gamma \in \dubbub_{k,s_{\indexc} - 16 - k}  \Big) \label{e.conindexbound}  \\
 & & \qquad  \, - \, \, \,  
\big( \conindex \log s_\indexc \big)^{-1} 4^{-1} s_\indexc \cdot  \conindex s_{\indexc}^{-2 \overline\theta - 11} \psap_{s_{\indexc}} \Big( \Gamma_{[0,r_j]} = \phi \Big)  \, . \nonumber 
\end{eqnarray}

Taking $j' = \indexc$ in~(\ref{e.keyob}) and using Proposition~\ref{p.rgjpnumber}, we learn that the summand in the first line satisfies 
\begin{eqnarray}
 &  &  \psap_{s_{\indexc}} \Big( \Gamma_{[0,r_j]} = \phi \, , \, \Gamma \in \dubbub_{k,s_{\indexc} - 16 - k} \Big) \label{e.summandeq} \\
 & = &  \psapleft{k} \Big( \Gamma_{[0,r_j]} = \phi \Big) \psap_{s_{\indexc}} \Big( \Gamma \in \dubbub_{k,s_{\indexc} - 16 - k} \Big) \nonumber \\
 & = & 
\psapleft{k} \Big( \Gamma_{[0,r_j]} = \phi \Big) \cdot p^{-1}_{s_{\indexc}} \cdot \lfloor \tfrac{1}{10} \, k^{1/2}  \rfloor \cdot \big\vert \sapell_k \big\vert \, \big\vert \sapr_{s_{\indexc} - 16 - k} \big\vert \, , \nonumber
\end{eqnarray} 
and thus by Lemma~\ref{l.polysetbound} (as well as $k \geq 2^{i+2}$ and $i \geq 7$) is at least
$$ 
\psapleft{k} \Big( \Gamma_{[0,r_j]} = \phi \Big) \cdot p^{-1}_{s_{\indexc}} \cdot \tfrac{1}{20} \, k^{1/2}  \cdot
 \tfrac{1}{8} p_k \cdot \tfrac{1}{2} p_{s_{\indexc} - 16 - k} \, . 
$$

Returning to~(\ref{e.conindexbound}), we find that
  $\psap_{s_{\indexc}} \big( \Gamma_{[0,r_j]} = \phi  \big)$ is at least
\begin{eqnarray}
 & & \Big( 1 +  
4^{-1}  \big( \log s_\indexc \big)^{-1}    s_{\indexc}^{-2 \overline\theta - 10} \Big)^{-1}  \big( \conindex \log s_\indexc \big)^{-1} \tfrac{1}{320} \, p^{-1}_{s_{\indexc}} \nonumber \\
& & \qquad \qquad \qquad \qquad \qquad \times \, \, \, \sum_{k \in \reg_{i+2}} \psapleft{k} \Big( \Gamma_{[0,r_j]} = \phi  \Big)
   k^{1/2} p_k p_{s_{\indexc} - 16 - k} \label{e.willbesame} \\
  & \geq &   \big( \conindex \log s_\indexc \big)^{-1} \tfrac{1}{400} \, p^{-1}_{s_{\indexc}} \mu^{s_\indexc - 16} \nonumber \\
& & \qquad \qquad \qquad \times \, \, \, \sum_{k \in \reg_{i+2}} \psapleft{k} \Big( \Gamma_{[0,r_j]} = \phi  \Big)
   k^{1/2} p_k \mu^{-k} (s_{\indexc} - 16 - k)^{-\thet - \eps} \, . \nonumber
\end{eqnarray}
Using $p_{s_{\indexc}} \leq \mu^{s_{\indexc}} s_{\indexc}^{-\thet + \eps}$, we obtain Lemma~\ref{l.compare}(1). 
   
\medskip

\noindent{{\em (2).}} For $\phi \in \fpart_{r_j}$,
\begin{eqnarray*}
 & &   \psap_{s_j} \Big(  \Gamma_{[0,r_j]} = \phi \, , \,  \Gamma \in \dubbub_{k,s_j - 16 - k} \, \, \textrm{for some $k \in \reg_{i+2}$} \Big) \\
 & \leq & \sum_{k \in \reg_{i+2}} \psap_{s_j}  \Big(  \Gamma_{[0,r_j]} = \phi \, , \,   \Gamma \in \dubbub_{k,s_j - 16 - k} \Big) \, .
 \end{eqnarray*}
By taking $j' = j$ in (\ref{e.keyob}), we see that the two equalities in the display~(\ref{e.summandeq}) hold with $s_j$ in place of $s_\indexc$. Since  $\big\vert \sapell_k \big\vert \leq p_k$ and $\big\vert \sapr_{s_j - 16 - k} \big\vert \leq p_{s_j - 16 - k}$, the last sum is seen to be at most
 \begin{eqnarray}
&  &  \tfrac{1}{10} \, p_{s_j}^{-1} \sum_{k \in \reg_{i+2}} 
   \psapleft{k}  \Big(  \Gamma_{[0,r_j]} = \phi \Big)  k^{1/2} p_k p_{s_j - 16 - k} \label{e.willbesametwo} \\
 & \leq &  \tfrac{1}{10} s_j^{\thet + \eps} \mu^{-16} \sum_{k \in \reg_{i+2}} 
   \psapleft{k}  \Big(  \Gamma_{[0,r_j]} = \phi \Big)  k^{1/2} p_k \mu^{-k} (s_j - 16 - k)^{-\thet + \eps} \nonumber \\
 & \leq &  \tfrac{1}{10}  \big( \tfrac{16}{3} \big)^{\thet} \mu^{-16}  \cdot s_j^{2\eps}  \sum_{k \in \reg_{i+2}} 
   \psapleft{k}  \Big(  \Gamma_{[0,r_j]} = \phi \Big)  k^{1/2} p_k  \mu^{-k}  \, , \nonumber   
 \end{eqnarray}
where we used $2^{i+5} \geq s_j \geq 2^{i+4}$, $k \leq 2^{i+3}$ and $i \geq 3$
in the form $s_j - 16 - k \geq 3s_j/16$.    \qed
  
\medskip

\noindent{\bf Proof of Proposition~\ref{p.similarity}.}
A consequence of Lemmas~\ref{l.dcub}(3) and~\ref{l.compare}. \qed

\subsection{Obtaining Proposition~\ref{p.almost} via similarity of measure}
This subsection is devoted to the next proof.
\medskip

\noindent{\bf Proof of Proposition~\ref{p.almost}.}
Note that the right-hand side of the bound
$$
  \psap_{s_{\indexc}} \Big( \Gamma_{[0,r_j]} \in \phicl{r_j}{s_j -1}{1/2 + 2\chi + \eclo} \Big) 
  \geq  
  \psap_{s_{\indexc}} \Big( \Gamma_{[0,r_j]} \in \phicl{r_j}{s_j -1}{1/2 + 2\chi + \eclo} 
  \, , \,  
\Gamma_{[0,r_j]} \in D_j  \Big)  
$$
is, by the similarity of measure Proposition~\ref{p.similarity}, at least the product of the quantity $c  \big( \log s_\indexc \big)^{-1} s_{\indexc}^{-4\eps}$ and the probability
$$
 \psap_{s_j} \Big( \Gamma_{[0,r_j]} \in \phicl{r_j}{s_j -1}{1/2 + 2\chi + \eclo} \, , \,  \Gamma \in \dubbub_{k,s_j - 16 - k} \, \, \textrm{for some $k \in \reg_{i+2}$} \, , \,  
\Gamma_{[0,r_j]} \in D_j   \Big) \, ,
$$
where the constant $c$ equals $\tfrac{1}{40}  \conindex^{-1} \big( \tfrac{3}{16} \big)^{\thet}$.
Note further that
\begin{itemize}
\item
$\psap_{s_j} \Big( \Gamma_{[0,r_j]} \in \phicl{r_j}{s_j -1}{1/2 + 2\chi + \eclo} \Big) \geq  1 - (s_j - 1)^{-(\chi + \eclo - \eroom)}$, 
by (\ref{e.acondition}) with $a=1$; and
\item
$\psap_{s_j} \Big(   \Gamma \notin \dubbub_{k,s_j - 16 - k} \, \, \textrm{for every $k \in \reg_{i+2}$}   \Big) \leq 1 - s_j^{-(\chi + \eclo/2 )}$, by Lemma~\ref{l.polyepsregthet} with $\varphi = \eclo/2$.
\end{itemize}
Also recall~(\ref{e.similarity}).
Thus,
\begin{eqnarray*}
 & & \psap_{s_{\indexc}} \Big( \Gamma_{[0,r_j]} \in \phicl{r_j}{s_j -1}{1/2 + 2\chi + \eclo} \Big) \\
 & \geq & c \, \big( \log s_\indexc \big)^{-1} s_{\indexc}^{-4\eps}  \, \bigg(  s_j^{-(\chi + \eclo/2 )}  - (s_j - 1)^{-(\chi + \eclo - \eroom)} \, - \, 4 \mu^{16} s_{\indexcj}^{-10}  \bigg) \, . 
\end{eqnarray*}

Recall that $s_j \geq 2^{i+4}$. Provided that the index $i \in \N$ is supposed to be sufficiently high, the displayed right-hand side is thus at least $\tfrac{1}{2} c  \big( \log s_\indexc \big)^{-1} s_{\indexc}^{-\chi - \eclo/2 - 4\eps}$ if we insist that $\eclo > 2\eroom$ as well as $\chi + \eclo < 10$ (the latter following from $\chi < 1/6$ and the harmlessly assumed $\eclo < 1$).

We are now ready to obtain Proposition~\ref{p.almost}.
We must set the values of the quantities $m$ and $m'$ in the proposition's statement. We take $m = s_{\indexc}$
and $m'$ equal to the constant difference $s_j - r_j$
between the terms in any pair in the sequence constructed in Definition~\ref{d.rjsj}. Recalling the bound on this difference stated in the definition, we see that 
 $2^{i+3} \leq m' \leq 2^{i+4} \leq m \leq 2^{i+5}$. 
 Set $\eclo = 4\eroom$ and $\eroom = \eps$, so that  $\chi + \eclo/2 + 4\eps = \chi + 6\eps < \esm - \eps$ and $1/2 + 2\chi + \eclo = \alpha$.  
 Note that $r_j$ belongs to the set $K$ specified in Proposition~\ref{p.almost} whenever $j \in [1,\indexc]$, because $r_j = s_j - m' \leq m - m'$. Thus, $\# K \geq 2^{i-4} \geq 2^{-9} m$, and 
Proposition~\ref{p.almost} is proved. \qed

\section{Step three: completing the proof of Proposition~\ref{p.thetexist}}\label{s.completing}
As we explained in this step's outline, we set the snake method index parameters
at the start of step three.  For a given sufficiently high  choice of $i \in \N$,  Proposition~\ref{p.almost} specifies  $m$ and $m'$. We then set the two index parameters, with $n$ equal to  $m-1$ and $\ell$ equal to $m - m'$.

For $\gamma \in \sap_{n+1}$, set 
$$
X_\gamma =  \sum_{k=0}^\ell  1\!\!1_{\gamma_{[0,k]} \in \phicl{k}{m'+k-1}{\alpha}} \, .
$$

Recall that $\gamma_{[0,\ell]} \in \fpart_{\ell,n}$ is $(\alpha,n,\ell)$-charming at an index $k \in [1,\ell]$ if
$$
 \psaw^0_{k + n - \ell} \Big( \Ga \closes \, \Big\vert \, \big\vert \Ga^1 \big\vert = k \, , \, \Ga^1 = \ga_{[0,k]} \Big) > n^{-\alpha} \, ;
$$
on the other hand, for given $k \in [1,\ell]$, the property that such a $\gamma$ satisfies $\ga_{[0,k]} \in \hPhi_{k,n-\ell+k}^\alpha$ takes the same form with $n^{-\alpha}$ replaced by the slightly larger quantity $(n - \ell + k)^{-\alpha}$. Since $n - \ell = m' - 1$ and $\beta = 1$, we find that, for any $\gamma \in \sap_{n+1}$,  
\begin{equation}\label{e.relate}
 X_\gamma \geq n^{\beta - \esm}/4 \, \, \, \textrm{implies that} \, \, \gamma_{[0,\ell]} \in \mathsf{CS}_{\beta,\esm}^{\alpha,\ell,n} \, .
\end{equation}

We now show that
\begin{equation}\label{e.claimx}
 \psap_{n+1} \Big( X_\Gamma \geq   n^{1- \esm + \eps/2}   \Big) \geq    n^{-\esm + \eps/2} \, .
\end{equation}

To derive this,  consider the expression
$$
S = \sum_{\gamma \in \sap_{n+1}}  \sum_{k=0}^{\ell}
1\!\!1_{\gamma_{[0,k]} \in \phicl{k}{m'+k-1}{\alpha}} \, .
$$

Recall that $p_{n+1}$ denotes $\# \sap_{n+1}$; using Proposition~\ref{p.almost} in light of~$n + 1 = m$,
\begin{eqnarray*}
  S  & = & 
 p_{n+1} \, \sum_{k=0}^{\ell}  \psap_{n+1} \Big(  \Gamma_{[0,k]} \in \phicl{k}{m'+k-1}{\alpha}   \Big) \\
 & \geq &   p_{n+1} \cdot \# K \cdot (n+1)^{-\esm + \eps} \geq  p_{n+1} \cdot 2^{-9} n \cdot  2^{-1} n^{-\esm + \eps}  \, .  
\end{eqnarray*}

Let $q$ denote the left-hand side of~(\ref{e.claimx}).
 Note that
$$
  S  \leq  p_{n+1}   \cdot \Big( q \big( \ell + 1 \big) + (1-q) n^{1-\esm + \eps/2}   \Big) \, .
$$
From the lower bound on $S$, and $n \geq \ell$,
$$
 q  \big( n + 1 \big) +  n^{1-\esm + \eps/2}
 \geq   2^{-10}  \, n^{1- \esm + \eps} \, ,
$$
which implies  for sufficiently high $n$ that $q \geq  n^{-\esm + \eps/2}$; in this way, we obtain~(\ref{e.claimx}).

Trivially, $n^{1-\esm + \eps/2} \geq n^{1 - \esm}/4$ (for all $n \in \N$, including our choice of $n$). Since the snake length exponent $\beta$ is set to one, we learn from~(\ref{e.relate}) and~(\ref{e.claimx}) that
$$
  \psap_{n+1} \Big( \Gamma_{[0,\ell]} \in \mathsf{CS}_{\beta,\esm}^{\alpha,\ell,n}   \Big) \geq \psap_{n+1} \Big( X_\Gamma \geq   n^{1-\esm + \eps/2}   \Big) \geq    n^{-\esm + \eps/2} \, .
 $$
Thus, if the dyadic scale parameter $i \in \N$ is chosen so that $n \geq 2^{i+4}$ is sufficiently high,  the charming snake presence hypothesis~(\ref{e.afewcs}) is satisfied. By Theorem~\ref{t.snakesecondstep}, $\psaw_n \big( \Gamma \closes \big) \leq 2(n+1)c^{-n^{\delta}/2}$. This deduction has been made for some value of $n \in (2\N + 1) \cap \big[2^{i+4}-1,2^{i+5}-1]$, where here $i \in \N$ is arbitrary to the right of a finite interval. Relabelling $c > 1$ to be any value in $(1,c^{1/2})$
completes the proof of Proposition~\ref{p.thetexist}. \qed

\section{Preparing for the proof of Theorem~\ref{t.thetexist}(1)}\label{s.prepare}

This argument certainly fits the template offered by the proof of Theorem~\ref{t.thetexist}(2). We will present the proof by explaining how to modify the statements and arguments used in the preceding derivation.

Recall first of all that the hypotheses of Theorem~\ref{t.thetexist}(2) implied the existence of the closing exponent, and that we chose to label $\chi \in \R$ such that $\psap_n\big( \Gamma \closes \big) = n^{-1/2 - \chi + o(1)}$ for odd~$n$;
 moreover, $\chi$ could be supposed to be non-negative. 
 Specifically, it was the lower bound  $\psap_n\big( \Gamma \closes \big) \geq n^{-1/2 - \chi - o(1)}$ that was invoked, in part one of the preceding proof.  

In the present proof, the closing exponent is not hypothesised to exist, and so we must abandon this usage of $\chi$. However, this parameter will be used again, with the basic role of~$1/2 + \chi$ as the exponent  in such a closing probability lower bound  being maintained.

\begin{definition}\label{d.highclose}
For $\maceta > 0$, define the set of indices $\highclose_\maceta \subseteq 2\N$ of $\maceta$-{\em high closing probability},
$$
 \highclose_\maceta = \Big\{ n \in 2\N : \psaw_{n-1} \big( \Gamma \closes \big) \geq n^{- \maceta} \Big\} \, .
$$
\end{definition}

For $\chi > 0$ arbitrary, we introduce the {\em closing probability}

\medskip

\noindent{\bf Hypothesis  $\mathsf{CP}_{\chi}$.} The set $2\N \setminus \hcb$ has limit supremum density in $2\N$ less than $1/1250$. 

\medskip

We begin the proof by the same type of reduction as was used in the earlier derivation. Here is the counterpart to  Proposition~\ref{p.thetexist}.
\begin{proposition}\label{p.closingprob}
Let $d = 2$. Let $\chi \in (0,1/{14})$, and assume Hypothesis~$\mathsf{CP}_\chi$. For some $c > 1$ and $\delta >0$, the set of odd integers $n$ satisfying $n + 1 \in \highclose_{1/2 + \chi}$ and 
$$
\psaw_n \big( \Ga \closes \big) \leq c^{-n^{\delta}}
$$
intersects the dyadic scale $\big[ 2^i,2^{i+1} \big]$ for all but finitely many $i \in \N$.
\end{proposition}
\noindent{\bf Proof of Theorem~\ref{t.thetexist}(1).}
The two properties of the index $n$ asserted by the proposition's conclusion are evidently in  contradiction for $n$ sufficiently high. Thus  Hypothesis~$\mathsf{CP}_\chi$ is false whenever  $\chi \in (0,1/{14})$. This implies the result. \qed

\medskip

The rest of the article is dedicated to proving Proposition~\ref{p.closingprob}. The three step plan of attack for deriving Proposition~\ref{p.thetexist} will also be adopted for the new proof.
Of course some changes are needed. 
We must cope with a deterioration in the known regularity of the $\theta$-sequence.
We begin by overviewing 
what we know in the present case and how the new information will cause changes in the three step plan.

In the proof of Theorem~\ref{t.thetexist}(2) via Proposition~\ref{p.thetexist}, we had the luxury of 
assuming the existence of the closing exponent. Throughout the proof of Proposition~\ref{p.closingprob} (and thus henceforth), we instead {\em fix $\chi \in (0,1/{14} )$, and suppose that Hypothesis~$\mathsf{CP}_\chi$ holds}. Thus, our given information is that the closing probability is at least $n^{-1/2 - \chi}$ for  a uniformly positive proportion of indices $n$ in any sufficiently long initial interval of positive integers.

We may state and prove right away the $\thet$-sequence regularity offered by this information.

\begin{proposition}\label{p.thetint}
Let $\eps \in (0,1)$. For all but finitely many values of $i \in \N$,
$$
 \# \Big\{ j \in 2\N \, \cap \, \big[2^{i},2^{i+1} \big]:  3/2  - \eps \leq \thet_ j \leq 3/2 + \chi + \eps \Big\}  \geq 2^{i-1} \big( 1 - \tfrac{1}{300} \big) \, .
$$
\end{proposition}
\noindent{\bf Proof.}
Theorem~\ref{t.polydev} implies that,  for all but finitely many values of $i \in \N$,
\begin{equation}\label{e.thirtyone}
 \# \Big\{ j \in 2\N \, \cap \, \big[2^{i},2^{i+1} \big]: \thet_ j \leq  3/2   - \eps \Big\} \leq   \tfrac{1}{600} \,  2^{i-1}  \, .
\end{equation}
Thus, it is enough to show that, also for all but finitely many $i$,
\begin{equation}\label{e.thirtytwo}
 \# \Big\{ j \in  2\N \, \cap \, \big[2^{i},2^{i+1} \big]: \thet_j > 3/2 + \chi + \eps \Big\} \leq   \tfrac{1}{600} \,  2^{i-1}  \, . 
\end{equation}

We will invoke Hypothesis $\mathsf{CP}_{\chi}$ as we verify~(\ref{e.thirtytwo}). 
First, we make a 

\medskip

\noindent{\bf Claim.}
For any
 $\varphi > 0$, the set  $\hcb  \cap   \big\{ n \in 2 \N: \thet_ n > 3/2 + \chi + \varphi \big\}$ is finite.

To see this, note that $2n p_n/c_{n-1} \geq n^{-1/2 - \chi}$ if $n \in \hcb$. From $c_{n-1} \geq \mu^{n-1}$ and $p_n = \mu^n n^{-\thet_n}$ follows $n^{-\thet_n} \geq \tfrac{1}{2\mu} n^{-3/2 - \chi}$ for such~$n$, and thus the claim.

The claim implies that, 
under Hypothesis $\mathsf{CP}_{\chi}$, the set
$\big\{ n \in 2\N : \thet_ n > 3/2 + \chi + \eps \big\}$
has limit supremum density  at most $1/1250$. From this,~(\ref{e.thirtytwo}) follows directly. \qed


\medskip

That is, the weaker regularity information forces a positive proportion of the polygon number deficit exponents $\thet_n$ to lie in any given open set containing the interval~$[3/2,3/2 + \chi]$. 
In the proof of Theorem~\ref{t.thetexist}(1), we knew this for the one-point 
set $\{ 3/2 + \chi \}$, for all high~$n$. 

In the three step plan, we will adjust step one so that the constant difference sequence constructed there incorporates the $\theta$-regularity that Proposition~\ref{p.thetint} shows to be typical, as well as a closing probability lower bound permitted by Hypothesis~$\mathsf{CP}_{\chi}$.

In step two, the weaker regularity will lead to a counterpart Lemma~\ref{l.polyepsreg} to the regulation polygon ample supply Lemma~\ref{l.polyepsregthet}. Where before we found a lower bound on the probability that a polygon is regulation global join of the form $n^{-\chi - o(1)}$, now we will find only a bound $n^{-2\chi - o(1)}$. 

Continuing this step, the mechanism of measure comparison for initial subpaths of polygons drawn from the laws $\psap_n$ for differing lengths~$n$ is no longer made via all regulation global join polygons but rather via such polygons whose length index lies in a certain regular set. The new similarity of measure result counterpart to Proposition~\ref{p.similarity} will be Proposition~\ref{p.similaritynew}. It will be applied at the end of step two to prove a slight variant of Proposition~\ref{p.almost}, namely Proposition~\ref{p.almostnew}.

This last result will yield Proposition~\ref{p.closingprob}
in step three by a verbatim argument to that by which Proposition~\ref{p.thetexist} followed from Proposition~\ref{p.almost}.

 \section{The three steps for Proposition~\ref{p.closingprob}'s proof in detail}\label{s.final}
 
In a first subsection of this section, we reset the snake method's parameters to handle the weaker information available.
In the following three subsections, we state and prove the assertions associated to each of the three steps.

\subsection{The snake method exponent parameters}\label{s.snakeparam}
Since $\chi < 1/{14}$, we may fix a parameter  $\eps \in \big(0,(1/2 - 7\chi)/14 \big)$, and do so henceforth. The three exponent parameters are then set so that
\begin{itemize}
\item $\beta = 1$;
\item $\alpha = 1/2 + 3\chi + 5\eps$; 
\item and $\esm = 4\chi + 9\eps$.
\end{itemize}
Note that the quantity $\delta = \beta - \esm - \alpha$, which must be positive if the method to work, is equal to $1/2 - 7\chi - 14\eps$.
The constraint imposed on $\eps$ ensures this positivity.

\subsection{Step one}

The essential conclusion of the first step one was the construction of the constant difference sequence in Definition~\ref{d.rjsj}. Now our aim is similar.

Fixing as previously two parameters $\eclo > \eroom > 0$ (which in the present case we will ultimately set $\eclo = 5\eroom$ and $\eroom = \eps$ in terms of our fixed parameter $\eps$), we adopt the Definition~\ref{d.cpt} of a closing probability typical index pair on dyadic scale~$i$, where of course the parameter $\chi$ is fixed in the way we have explained. An index pair $(k,j)$ that satisfies this definition is called {\em regularity typical} for dyadic scale~$i$ if moreover  $j \in \highclose_{1/2 + \chi}$ and  
\begin{equation}\label{e.tineqs}
 3/2  - \epn  \leq \thet_j \leq 3/2 + \chi + \epn  \, .
\end{equation}

We then replace Definition~\ref{d.rjsj} with the following.
\begin{definition}\label{d.rjsj.new}
Set $\indexc = 2^{i-4}$. We specify a {\em constant difference} sequence of  index pairs $(r_j,s_j)$, $1 \leq j \leq \indexc$, whose elements are regularity typical for dyadic scale $i$, with $s_j - r_j$ independent of $j$ and valued in $[2^{i+4}- 2^i,2^{i+4}]$,  and $\big\{ r_j: 1 \leq j \leq \indexc \big\}$ increasing. 
Moreover, $s_j \in 2\N \cap [ 2^{i+4} , 2^{i+4} + 2^i + 2^{i-2}  ]$ for each $j \in [1,\indexc]$.
\end{definition}
Our job in step one is to construct such a sequence.
First, recall the set $E$ of index pairs $(k,j)$ specified in Definition~\ref{d.e}. Let $E'$ be the subset of $E$ consisting of such pairs for which~$j \in \highclose_{1/2 + \chi}$ and~(\ref{e.tineqs}) holds.
Counterpart to Lemma~\ref{l.etypical}, we have:
\begin{lemma}\label{l.eprimetypical}
Provided that the index $i \in \N$ is sufficiently high, except
for at most $\tfrac{1}{100} 2^{i+3}$ values of $j \in  2 \N \cap \big[ 2^{i+4},  2^{i+5} \big]$, the set of $k \in \big[ 2^i , 2^i + 2^{i-2} \big]$ such that $(k,j) \not\in E'$ has cardinality at most $2 \big(  2^{i+5} \big)^{1-\eroom}$.
\end{lemma}
\noindent{\bf Proof.}
Proposition~\ref{p.thetint} implies that the set of $j \in 2\N \cap \big[ 2^{i+4},2^{i+5} \big]$ satisfying~(\ref{e.tineqs}) -- a set containing all second coordinates of pairs in $E'$ -- has cardinality at least $\big( 1 -\tfrac{1}{300} \big) 2^{i+3}$. Invoking Hypothesis~$\mathsf{CP}_\chi$, we see that, of these values of~$j$, at most~$(2^{i+4} + 1) \tfrac{1}{1250} \leq  \tfrac{1}{600} 2^{i+3}$ fail the test of membership of~$\highclose_{1/2 + \chi}$.
We claim that the set of remaining values of $j$ 
satisfies the statement in the lemma. Indeed, all these values
belong to $\highclose_{1/2 + \chi}$. As such, we are able to invoke Lemma~\ref{l.closecard}, as we did in proving Lemma~\ref{l.etypical}, but this time with $\alphamac = 1/2 + \chi$ and $\deltamac = \eroom$, to find that at most  $2 \big(  2^{i+5} \big)^{1-\eroom}$ values of  $k \in \big[ 2^i , 2^i + 2^{i-2} \big]$ violate the condition~(\ref{e.rchithet}). \qed

Analogously to Proposition~\ref{p.manyethet},
\begin{proposition}\label{p.manyethet.new}
There exists $k \in 2\N \cap \big[ 2^{i+4}, 2^{i+4} + 2^i \big]$ for which there are at least $2^{i-4}$ values of $j \in 2\N \cap [0,2^{i-2}]$ with the index pair
$\big( 2^i + j, k + j \big)$ regularity typical on dyadic scale $i$.
\end{proposition}
\noindent{\bf Proof.} Using Lemma~\ref{l.eprimetypical},
we find, analogously to~(\ref{e.ineqe}),
\begin{eqnarray*}
 & & \sum_{j \in 2\N \cap [ 2^{i+4} + 2^{i-2} , 2^{i+4} + 2^{i} ] } \sum_{\ell = 2^i}^{2^i + 2^{i-2}} 
 1\!\!1_{(\ell,j) \in E'} \\ 
 & \geq &  \Big( 2^{i-2} + 1 - 2 \big( 2^{i+5} \big)^{1 - \eroom} \Big) \cdot  \Big( \tfrac{1}{2} \big( 2^{i} - 2^{i-2} \big) - \tfrac{1}{100} 2^{i+3} \Big) \, . 
\end{eqnarray*}
The argument for Proposition~\ref{p.manyethet} now yields the result. \qed

\medskip

We may now conclude step one, because the proof of Proposition~\ref{p.notquite} is valid with the notion of regularity typical replacing closing probability typical; thus, we have specified a constant difference sequence in the sense of Definition~\ref{d.rjsj.new}. (Incidentally, the unused Proposition~\ref{p.notquite} is seen to be valid with right-hand side strengthened to $1 - 2^{-2i\chi}$ if we choose $\eclo = 5\eps$ and $\eclo > \eroom$ by considering $a=2$ in~(\ref{e.acondition}).)

\subsection{Step two}

The new version of Proposition~\ref{p.almost} differs from the original only in asserting that the constructed $m$ belongs to $\highclose_{1/2 + \chi}$.
\begin{proposition}\label{p.almostnew}
For each $i \in \N$ sufficiently high, there exist $m \in 2\N \cap \big[ 2^{i+4}, 2^{i+5} \big] \cap \highclose_{1/2 + \chi}$ and $m' \in \big[ 2^{i+3},2^{i+4} \big]$ such that, writing $K$ for the set of values $k \in \N$, $1 \leq k \leq m-m'$, that satisfy
$$
\psap_m \Big( \Ga_{[0,k]} \in \hPhi_{k,k+m'-1}^\alpha \Big)
 \geq  m^{-\esm + \eps} \, ,
$$
we have that $\vert K \vert \geq 2^{-9} m$.
\end{proposition}

In this step, we provide a suitable similarity of measure result, Proposition~\ref{p.similaritynew}, and use it to derive Proposition~\ref{p.almostnew}.

Before doing so, we must adapt our regulation polygon tools to cope with the weaker regularity that we are hypothesising. The ample regulation polygon supply Lemma~\ref{l.polyepsregthet}
will be replaced by Lemma~\ref{l.polyepsreg}. 

In preparing to state the new lemma, recall that Lemma~\ref{l.polyepsregthet} 
was stated in terms of a regular index set $\reg_{i+2}$ that was specified in a trivial way.  
Recalling  that $\indexc = 2^{i-4}$, we now specify a non-trivial counterpart set $\reg_{n,i+2}^{\chi + \eps}$ when $n$ is an element in the sequence $\big\{ s_j: 1 \leq j \leq \indexc \big\}$ from step one. 
\begin{definition}
For $j \in \big[1, \indexc \big]$,   define
$\regnew{s_j}{i+2}{\chi + \epn}$ to be the set of $k \in 2\N \cap \big[ 2^{i+2} , 2^{i+3} \big]$ such that 
$$
  \max \big\{ \thet_k , \thet_{s_j-16-k} ,  
\thet_{s_\indexc - 16 - k}  \big\}  \leq  3/2 +  \chi + \eps
\, \, \, \textrm{and} \, \, \,
\thet_{s_j - 16 - k} \geq 3/2 - \epn  \, \, .   
$$
\end{definition}

Here is the new ample supply lemma.
\begin{lemma}\label{l.polyepsreg}
Let $j \in [1,\indexc]$. 
Then
$$
\psap_{s_j} \Big( \Gamma \in \dubbub_{k,s_j-16-k} \, \, \textrm{for some $k \in  \regnew{s_j}{i+2}{\chi + \epn}$} \Big) \geq \frac{c}{\log s_j} \, {s_j}^{-(2\chi + 3\epn)} \, ,
$$
where $c > 0$ is a universal constant.
\end{lemma}

The similarity of measure Proposition~\ref{p.similarity} is replaced by the next result. 
\begin{proposition}\label{p.similaritynew}
For each $j \in [1,L]$, there is a subset $D_j \subseteq \fpart_{r_j}$ 
satisfying
\begin{equation}\label{e.similaritynew}
 \psap_{s_j} \Big( \Gamma_{[0,r_j]} \not\in D_j \, , \, \Gamma \in \dubbub_{k,s_j - 16 - k} \, \, \, \textrm{for some $k \in \regnew{s_j}{i+2}{\chi + \epnmac}$} \Big) \leq 4 \mu^{16}  s_j^{-10} 
\end{equation}
such that, for each $\phi \in D_j$,
\begin{eqnarray*}
  \psap_{s_{\indexc}} \Big( \Gamma_{[0,r_j]} = \phi  \Big) 
  & \geq  & \tfrac{1}{40} \,  \conindex^{-1} 
\big( \tfrac{3}{16} \big)^{3/2}  \big( \log s_\indexc \big)^{-1}  s_\indexc^{-2\chi -  4\epn} \\ 
 & &    \times \, \, \, \psap_{s_j} \Big(  \Gamma_{[0,r_j]} = \phi \, , \,  \Gamma \in \dubbub_{k,s_j - 16 - k} \, \, \textrm{for some $k \in \regnew{s_j}{i+2}{\chi + \epnmac}$} \Big) \, .
\end{eqnarray*}
\end{proposition}

In the rest of this subsection of step two proofs, we prove in turn Lemma~\ref{l.polyepsreg} and Propositions~\ref{p.similaritynew} and~\ref{p.almostnew}.

\subsubsection{Deriving Lemma~\ref{l.polyepsreg}}

We will use the next result.

\begin{lemma}\label{l.rsize}
For  $i \in \N$ with $i \geq 6$, and $n \in 2\N \cap \big[ 2^{i+4},2^{i+5} \big]$,
$$
  \big\vert \regnew{s_j}{i+2}{\chi + \epn}  \big\vert \geq  \big( \tfrac{9}{10} - \tfrac{1}{32} \big) 2^{i+1}    \, .
$$

\end{lemma}
\noindent{\bf Proof.} It is enough to prove two claims.

For $\varphi > 0$, set
$\regul_{s_j,i+2}^\varphi =     \big\{ k \in 2\N \, \cap \, \big[2^{i+2},2^{i+3} \big]:  \max \big\{ \thet_k , \thet_{s_j-16-k}  \big\}  \leq  3/2 +  \varphi \,  \big\}$. 

\noindent{\em Claim $1$.}
The set  $\regul_{s_j,i+2}^{\chi + \epn} \setminus \regnew{s_j}{i+2}{\chi + \epn}$ has cardinality at most $2^{i-4}$.

\noindent{\em Claim $2$.}
$\big\vert \regul_{s_j,i+2}^{\chi + \epn} \big\vert \geq \tfrac{9}{10} \, 2^{i+1}$.

\noindent{\em Proof of Claim~$1$.} Note that $k \in \regul_{s_j,i+2}^{\chi + \epn} \setminus \regnew{s_j}{i+2}{\chi + \epn}$ implies that either $s_j - 16 - k$ belongs to the union up to index~$i+4$ (in the role of~$i$) of the sets in~(\ref{e.thirtyone}), or $s_\indexc - 16 - k$ belongs to the comparable union of the sets in $(\ref{e.thirtytwo})$. (The reason that the index here is $i+4$ is because $s_j$ and $s_\indexc$ are at most $2^{i+4} + 2^i + 2^{i-2} \leq 2^{i+5}$ in view of Definition~\ref{d.rjsj.new}.) Thus, the cardinality in question is at most $\tfrac{4}{600} 2^{i+3}$, where the factor of four accomodates the two alternatives as well as the sum over indices of dyadic scales up to index~$i+4$.  
We obtain Claim~$1$.

\noindent{\em Proof of Claim~$2$.}
For $\varphi > 0$, define the set of $\varphi$-high $\thet$ values on dyadic scale $i$,
$$
\macess_{i}^\varphi = \Big\{ k \in 2\N \, \cap \, \big[2^{i},2^{i+1} \big]: \thet_k \geq  3/2  + \varphi \Big\}  \, ,
$$
and also set
$\irr_{s_j,i+2}^\varphi = 2\N \cap \big[ 2^{i+2} , 2^{i+3} \big] \setminus \regul_{s_j,i+2}^\varphi$.
We have then that,
 for any $\varphi > 0$,
$$
 \irr_{s_j,i+2}^\varphi \subseteq  \macess_{i+2}^\varphi \, \cup \, \Big\{ k \in 2\N \cap \big[ 2^{i+2} , 2^{i+3} \big]: \thet_ {s_j-16-k} \geq   3/2 + \varphi \Big\} \, ;
 $$
 the latter set in the right-hand union is a subset of
 $$ 
 2\N \cap \big[ 2^{i+3} - 2^4,2^{i+3} \big] \, \, \cup \, \, \Big\{ k \in 2\N \cap \big[ 2^{i+2},2^{i+3} - 2^4 \big]: s_j-16-k \in \macess_{i+3}^\varphi \cup \, \macess_{i+4}^\varphi \Big\} \, ,
$$
because $2^{i+3} \leq s_j-16-k \leq  2^{i+5}$ for $k \in \big[ 2^{i+2},2^{i+3} - 2^4 \big]$. Thus, 
$$
\# \irr_{s_j,i+2}^{\chi + \epn} \leq  \# \macess_{i+2}^{\chi + \epn} +  \# \macess_{i+3}^{\chi + \epn} + \# \macess_{i+4}^{\chi + \epn} + 9 \, ,
$$
 so that~(\ref{e.thirtytwo}) implies $\big\vert \irr_{s_j,i+2}^{\chi + \epn} \big\vert \leq \tfrac{7}{600} \, 2^{i+1} + 9 \leq \tfrac{1}{10}  \, 2^{i+1}$ (for $i \geq 6$). Thus, we obtain Claim~$2$. \qed

\medskip

\noindent{\bf Proof of Lemma~\ref{l.polyepsreg}.}
 The proof is similar to that of Lemma~\ref{l.polyepsregthet}. Now, the probability in question equals the ratio of $\big\vert \bigcup_{l \in  \regnew{s_j}{i+2}{\chi + \epn}} \dubbub_{l,s_j-16-l} \big\vert$ and $p_{s_j}$.

To bound the numerator below, we again seek to  apply Proposition~\ref{p.polyjoin.aboveonehalf} with the role of $i \in \N$ played by $i+2$, but now with $n = s_j$ and $\mathsf{R} =  \regnew{s_j}{i+2}{\chi + \epn}$. Set $\Theta$ in the proposition to be $3/2 + \chi + \eps$. Any element of $\regnew{s_j}{i+2}{\chi + \epn}$ may act as $\kmac \in \mathsf{R}$
in the way that the proposition demands, by the definition of $\regnew{s_j}{i+2}{\chi + \epn}$; (and Lemma~\ref{l.rsize} shows that $\regnew{s_j}{i+2}{\chi + \epn}$ is non-empty if $i \geq 6$).
Applying the proposition and then using Lemma~\ref{l.rsize}, $2^{i+1} \geq 2^{-4} s_j$,
and the definition of  $\regnew{s_j}{i+2}{\chi + \epn}$, we find that 
$$
 \bigg\vert \, \bigcup_{l \in  \regnew{s_j}{i+2}{\chi + \epn}} \dubbub_{l,s_j-16-l} \, \bigg\vert \geq c \, \frac{s_j^{1/2}}{\log s_j} \cdot  \big( \tfrac{9}{10} - \tfrac{1}{32} \big) 2^{-4} s_j  \cdot \mu^{s_j - 16} s_j^{-(3 + 2\chi + 2\eps)}  \, ,   
$$
which is at least a small constant multiple of  $\mu^{s_j} \big(\log s_j\big)^{-1}s_j^{-3/2 - 2\chi - 2\eps}$.

The ratio's denominator is bounded above with the bound  $p_{s_j} \leq \mu^{s_j} \, {s_j}^{-3/2  + \epn}$ which is due to the lower bound in~(\ref{e.tineqs}) with $s_j$ playing the role of $j$. The lemma follows by relabelling~$c > 0$. \qed

\medskip

\subsubsection{Deriving Proposition~\ref{p.similaritynew}.} We begin by respecifying the set $\tau_\gamma$.
For each $\gamma \in \sap_{s_j}$, $1 \leq j \leq \indexc$, we now set
$$
\tau_\gamma = \Big\{ k \in \regnew{s_j}{i+2}{\chi + \epnmac} : \gamma \in \dubbub_{k,s_{\indexcj} - 16 - k} \Big\} \, .
$$

Note that $\regnew{s_j}{i+2}{\chi + \epnmac}$ has replaced $\reg_{i+2}$ here. As we comment on the needed changes, it is understood that this replacement is always made.

We reinterpret the quantity $\overline{\theta}$ to be $3/2 + \chi + \epn$. The constant $\conindex > 0$ may again be chosen in order that~(\ref{e.constl}) be satisfied.
Indeed, since $(r_\indexc,s_\indexc)$ is by construction regularity typical on dyadic scale~$i$, we have that $s_\indexc$  is an element of $\mathsf{HCP}_{1/2 + \chi}$ and thus also of $\mathsf{HPN}_{3/2 + \chi + o(1)}$.
 Corollary~\ref{c.globaljoin} is thus applicable with $n = s_\indexc$, so that~(\ref{e.constl}) results.

The set $D_j$ for $j \in [1,\indexc]$ is specified as before, and Lemma~\ref{l.dcub}(1) is also obtained, by the verbatim argument. 

The definition of $\overline\theta$ has been altered in order to maintain the validity of Lemma~\ref{l.dcub}(3). The moment in the proof of this result where the new definition is needed is in deriving the bound 
\begin{equation}\label{e.boundadapt}
p_{s_\indexc}^{-1} p_k p_{s_\indexc - 16 - k} \geq \mu^{-16} s_\indexc^{-2\overline{\theta}}
\end{equation}
for $k \in \regnew{s_j}{i+2}{\chi + \epnmac}$
 in the proof of  Lemma~\ref{l.dcub}(2) (which is needed to prove the third part). That  $k \in \regnew{s_j}{i+2}{\chi + \epnmac}$ ensures that both $\thet_k$ and $\thet_{s_\indexc - 16 - k}$ are at most $3/2 + \chi + \epnmac = \overline{\theta}$, from which follows~(\ref{e.boundadapt}), and thus Lemma~\ref{l.dcub}(2) and~(3).
 
The weaker form of the inference that we may make in view of the weaker available regularity information is evident in the counterpart to Lemma~\ref{l.compare}. Note the presence of the new factors of $s_\indexc^{-\chi}$ and $s_j^{\chi}$ in the two right-hand sides.
\begin{lemma}\label{l.comparenew}
Recall that $\indexc = 2^{i-4}$.  Let $j \in [1,\indexc]$.
\begin{enumerate}
\item 
For $\phi \in D_j$,
\begin{eqnarray*}
   \psap_{s_{\indexc}} \Big( \Gamma_{[0,r_j]} = \phi  \Big) 
 & \geq &
      \tfrac{1}{400} \, \conindex^{-1} \big( \log s_\indexc \big)^{-1}  s_{\indexc}^{-\chi - 2\epn} \mu^{-16} \\
      & & \qquad \qquad \times \, \, \, \sum_{k \in  \regnew{s_j}{i+2}{\chi + \epn}} \psapleft{k} \Big( \Gamma_{[0,r_j]} = \phi  \Big)
   k^{1/2} p_k \mu^{-k} \, .
\end{eqnarray*}
\item  For $\phi \in \fpart_{r_j}$,
\begin{eqnarray*}
   & & \psap_{s_j} \Big(  \Gamma_{[0,r_j]} = \phi \, , \,  \Gamma \in \dubbub_{k,s_j - 16 - k} \, \, \textrm{for some $k \in  \regnew{s_j}{i+2}{\chi + \epn}$} \Big) \\
   & \leq &  \tfrac{1}{10}  \big( \tfrac{16}{3} \big)^{3/2} \mu^{- 16} \, \cdot s_j^{\chi + 2\epn}  \sum_{k \in \regnew{s_j}{i+2}{\chi + \epnmac}} 
   \psapleft{k}  \Big(  \Gamma_{[0,r_j]} = \phi \Big)  k^{1/2} p_k \mu^{-k} \, .  
\end{eqnarray*}
\end{enumerate}
\end{lemma}
\noindent{\bf Proof (1).} 
The earlier proof runs its course undisturbed until 
$\psap_{s_{\indexc}} \big( \Gamma_{[0,r_j]} = \phi  \big)$
is found to be at least the expression in the double line ending at~(\ref{e.willbesame}). The proof of this part is then completed by noting the bound
$$
 p_{s_\indexc}^{-1} p_{s_\indexc - 16 - k} \geq \mu^{-16 - k} s_\indexc^{-\chi - 2\epn} 
$$
for $k \in  \regnew{s_j}{i+2}{\chi + \epn}$.
This bound is due to $\thet_{s_\indexc - 16 -k} \leq 3/2 + \chi + \epn$ and $\thet_{s_\indexc} \geq 3/2  - \epn$, the latter being the lower bound in~(\ref{e.tineqs}) with $j = s_\indexc$.

\medskip

\noindent{{\bf (2).}} 
The left-hand side of the claimed inequality is bounded above by~(\ref{e.willbesametwo}). The proof is then completed by invoking the bound
\begin{equation}\label{e.neededbound}
 p_{s_j}^{-1} p_{s_j - 16 - k} \leq \mu^{-16 - k} \big( \tfrac{16}{3} \big)^{3/2} s_j^{\chi + 2\epn}
\end{equation}
for $k \in  \regnew{s_j}{i+2}{\chi + \epn}$.
To derive this bound, note that
$$
 p_{s_j - 16 -k} \leq \mu^{s_j - 16 - k} \big( s_j - 16 - k \big)^{-3/2  + \epn} \leq \mu^{s_j - 16 -k} \big( \tfrac{16}{3} \big)^{3/2} s_j^{-3/2  + \epn} \, , 
$$
where the first inequality follows from 
$\thet_{s_j - 16 -k} \geq 3/2  - \epn$, which is due to $k \in   \regnew{s_j}{i+2}{\chi + \epn}$;
and where the second depends on  $s_j - 16 - k \geq 3s_j/16$, which we saw at the corresponding moment in the derivation that is being adapted, as well as on $0 \leq \epn \leq 3/2$.

We also have  
$$
 p_{s_j} \geq \mu^{s_j} s_j^{-3/2 - \chi - \epn}
$$
due to the upper bound in~(\ref{e.tineqs}) with the role of $j$ played by $s_j$.

The last two displayed inequalities combine to yield~(\ref{e.neededbound}). \qed

\medskip

\noindent{\bf Proof of Proposition~\ref{p.similaritynew}.}
A consequence of Lemmas~\ref{l.dcub}(3) and~\ref{l.comparenew}. \qed

\subsubsection{Proof of Proposition~\ref{p.almostnew}.}
Note that the right-hand side of the bound
$$
  \psap_{s_{\indexc}} \Big( \Gamma_{[0,r_j]} \in \phicl{r_j}{s_j -1}{1/2 + 3\chi + \eclo} \Big) 
  \geq  
  \psap_{s_{\indexc}} \Big( \Gamma_{[0,r_j]} \in \phicl{r_j}{s_j -1}{1/2 + 3\chi + \eclo} 
  \, , \,  
\Gamma_{[0,r_j]} \in D_j  \Big)  
$$
is, by Proposition~\ref{p.similaritynew}, at least the product of $c_1  \big( \log s_\indexc \big)^{-1} s_\indexc^{-2\chi -  4\epn }$ and
\begin{equation}\label{e.hfrg}
 \psap_{s_j} \Big( \Gamma_{[0,r_j]} \in \phicl{r_j}{s_j -1}{1/2 + 3\chi + \eclo} \, , 
    \Gamma \in \dubbub_{k,s_j - 16 - k} \, \, \textrm{for some $k \in \regnew{s_j}{i+2}{\chi + \epn}$} \, , \,  
\Gamma_{[0,r_j]} \in D_j   \Big)
\end{equation}
where $c_1 = \tfrac{1}{40} \,  \conindex^{-1}  \big( \tfrac{3}{16} \big)^{3/2}$.

Recall that the index pair $(r_j,s_j)$ is regularity, and thus closing probability, typical on dyadic scale~$i$. We may thus apply 
(\ref{e.acondition}) with $s_j$ playing the role of $j$ and with $a=2$ to find that 
\begin{itemize}
\item
$\psap_{s_j} \Big( \Gamma_{[0,r_j]} \in \phicl{r_j}{s_j - 1}{1/2 + 3\chi + \eclo} \Big) \geq  1 - (s_j - 1)^{-(2\chi + \eclo - \eroom)}$.
\end{itemize}

Also recall 
\begin{itemize}
\item 
that, by  Lemma~\ref{l.polyepsreg}, the quantity 
$\psap_{s_j} \Big(   \Gamma \notin \dubbub_{k,s_j - 16 - k} \, \, \textrm{for every $k \in \regnew{s_j}{i+2}{\chi + \epn}$}   \Big)$ is at most $1 \, - \, c \, \big( \log s_j \big)^{-1}  s_j^{-(2\chi + 3\epn)}$; and
\item the inequality (\ref{e.similaritynew}).
\end{itemize}

The three inequalities stated or recalled in these bullet points may be used to find a lower bound on the expression~(\ref{e.hfrg}).
 We see then that
\begin{eqnarray*}
 & & \psap_{s_{\indexc}} \Big( \Gamma_{[0,r_j]} \in \phicl{r_j}{s_j - 1}{1/2 + 3\chi + \eclo} \Big) \\
 & \geq & c_1  \big( \log s_\indexc \big)^{-1} s_{\indexc}^{-2\chi -  4\epn}  \,  
   \, \bigg(  \, \frac{c}{\log s_j} \, s_j^{-(2\chi + 3\epn)} \, - \, (s_j - 1)^{-(2\chi + \eclo - \eroom)} \, - \,  4 \mu^{16} s_{\indexcj}^{-10} \, \bigg) \, . 
\end{eqnarray*}

Provided that the index $i \in \N$ is supposed to be sufficiently high, the right-hand side is at least $\tfrac{1}{2} c_1 c \big( \log s_\indexc \big)^{-2} s_{\indexc}^{-4\chi - 7\eps} \geq s_{\indexc}^{-4\chi - 8\eps} = s_{\indexc}^{-\eta + \eps}$ if we insist that $\eclo - \eroom > 3 \epn$ as well as $2\chi + \eclo - \eroom < 10$.
We may set $\eclo = 5\eroom$ and $\eroom = \eps$, because we imposed at the outset that $\chi \in (0,1/{14})$ and in Section~\ref{s.snakeparam} the condition  $\epn \in \big(0,(1/2-7\chi)/14\big)$, and these choices do ensure the stated bounds. Note that $\alpha = 1/2 + 3\chi + \eclo$ results from this choice.

We may now obtain  
Proposition~\ref{p.almostnew}. We set the values of $m$ and $m'$ as we did in proving Proposition~\ref{p.almost}, taking  $m = s_{\indexc}$ and 
$m'$ equal to the constant difference $s_j - r_j$ associated to any element in the sequence $(r_j,s_j)$ from Definition~\ref{d.rjsj.new} (in place of course of Definition~\ref{d.rjsj}).
The bound $\# K \geq 2^{i-4} \geq 2^{-9} m$ holds as in the earlier proof. Note however that we are also claiming that  $m \in \highclose_{1/2 + \chi}$.  
This holds because $m = s_\indexc$ is the second coordinate of a regularity typical index pair for dyadic scale~$i$.  
Proposition~\ref{p.almostnew} is proved. \qed

\subsection{Step three}

At the start of the earlier step three, the index parameters $n$ and $\ell$ were respectively set equal to $m-1$ and $m -m'$, where Proposition~\ref{p.almost} provided $m$ and $m'$. We now make the same choice, with Proposition~\ref{p.almostnew} providing the latter two quantities.

\medskip

\noindent{\bf Proof of Proposition~\ref{p.closingprob}.} The completion of the proof of Proposition~\ref{p.thetexist} in Section~\ref{s.completing} applies verbatim after recalling that the snake method's parameter $\delta = \beta - \esm - \alpha$ is positive, and noting that the method's index parameter $n$ satisfies $n+1 \in \highclose_{1/2 + \chi}$
because it has been set equal to~$m - 1$ where $m$ is specified in Proposition~\ref{p.almostnew}. \qed

\bibliographystyle{plain}

\bibliography{saw}

\end{document}